\documentclass[10pt]{article}
\usepackage{latexsym}
\usepackage{amsfonts}
\usepackage{enumitem}
\usepackage[active]{srcltx} 
\usepackage{color}

\usepackage{amsmath}
\usepackage{amssymb}

\usepackage{slashbox}

\newcommand\tstrut{\rule{0pt}{3.2ex}}

\newcommand{\ba}{\begin{array}}\newcommand{\ea}{\end{array}}
\newcommand{\sym}{S}
\newcommand{\ns}{\rm}

\newcommand{\nse}{\kern-3pt\ns$=$}\newcommand{\qd}{\hfill$\Box$\medbreak}

\newcommand{\noi}{\noindent}
\newcommand{\ext}{\raise1pt\hbox{$\ts\bigwedge$}}
\renewcommand{\sym}{S}
\newcommand{\ts}{\textstyle}
\newcommand{\rf}[1]{(\ref{#1})}
\newcommand{\chii}{\raise2pt\hbox{$\chi$}}

\newcommand{\Fg}{\mbox{${\cal F}\kern-2pt_g$}}

\newcommand{\Mg}{\mbox{${\cal M}\kern-2pt_g$}}
\newcommand{\Ng}{\mbox{${\cal N}\kern-2pt_g$}}
\newcommand{\V}{V\kern-1pt}

\newcommand{\Gg}{\mbox{${\cal G}\kern-2pt_g$}}
\newcommand{\cir}{\raise1.6pt\hbox{\footnotesize$\circ$}}

\newcommand{\Res}[2]{\hbox{\ns Res}\kern-16pt\lower5pt\hbox{\footnotesize$_{#1}$}\kern2pt\left[#2\right]}
\newcommand{\qk}{quaternion-K\"ahler\kern2pt}\renewcommand{\,}{\kern1pt}

\newcommand{\dirac}{/\kern-5pt\partial}

\newcommand{\End}{{\rm End}}

\newcommand{\lra}{\longrightarrow}
\renewcommand{\ts}{\textstyle}
0
\newtheorem{theo}{Theorem}[section]
\newtheorem{defi}{Definition}[section]

\def\frac#1#2{{#1\over#2}}
\def\be#1\ee{\begin{equation}#1\end{equation}}

\begin{document}
\title{A note on the geometry and topology of almost even-Clifford Hermitian manifolds}

\author{
Gerardo Arizmendi\footnote{Centro de
Investigaci\'on en Matem\'aticas, A. P. 402,
Guanajuato, Gto., C.P. 36000, M\'exico. E-mail: gerardo@cimat.mx, lucia@cimat.mx, rherrera@cimat.mx},
\,\,\,\,\,
Ana L. Garcia-Pulido$^*$
\,\,\,\,\, and
Rafael Herrera$^*$\footnote{Partially supported by a CONACYT grant
}
}

\date{}

\maketitle

{
\abstract{

We compute the structure groups of almost even-Clifford Hermitian manifolds and determine when such groups lead to Spin structures.

}
}

\section{Introduction}

Almost even-Clifford Hermitian structures on oriented Riemannian manifolds were introduced recently in their current form in \cite{Moroianu-Semmelmann} under the simpler name of
even Clifford structures. They are a subject of current interest \cite{Arizmendi-Hadfield, Arizmendi-Herrera-Santana, Dearricott,
Moroianu-Pilca, Piccinni}, although  similar types of structures have been studied in the past \cite{Barberis-Dotti-Miatello, Joyce,Nikolayevsky}.
They are generalizations of almost Hermitian and almost quaternion-Hermitian structures, and there has been quite some interest in them .
The existence of such a structure on a manifold 
implies the reduction of its structure group to the normalizer of the homomorphic image of a Spin group. In this paper, we identify such structure group  
(cf. Theorem \ref{theo: quotient groups}) by using the results about their Lie algebras given in \cite{Arizmendi-Herrera}. 
In the case of $4m$-dimensional almost quaternion-Hermitian manifolds, we
know that such manifolds are Spin when $m$ is even. This is due to the (topological) reduction of the structure group from $SO(4n)$ to 
the Lie group $Sp(n)Sp(1)$ which, in turn, embeds into $Spin(4m)$ when $m$ is even.
Thus, by analogy, we were led to study when such manifolds admit Spin structures.

Recall that
an oriented $n$-dimensional Riemannian manifold is Spin if its orthonormal frame bundle $P_{SO}$ admits
a double cover by a principal $Spin(n)$ bundle $P_{Spin}$ 
\[\Lambda:P_{Spin}\longrightarrow P_{SO}\]
which is $Spin(n)$ equivariant, i.e. $\Lambda(pg) = \Lambda(p) \lambda_n(g)$ for all $g\in Spin(n)$.
A Riemannian manifold will automatically be Spin if its structure group reduces to a proper  
subgroup $G\subset SO(n)$ such that there exists a lifting map which makes the following diagram commute
\[
\begin{array}{lll}
 &  & Spin(n)\\
 & \nearrow & \downarrow\\
G& \hookrightarrow & SO(n).
\end{array}
\]
Indeed, such a lift exists if and only if $\pi_1(G)$ maps trivially into $\pi_1(SO(n))$.

In this paper, we determine when there exists a lifting map which makes the following diagram commute (cf. Theorem \ref{theo: lifts})
\[
\begin{array}{lll}
 &  & Spin(N)\\
 &  & \downarrow\\
\mathcal{N}_{SO(N)}(S) & \hookrightarrow & SO(N).
\end{array}
\]
where $N$ stands for the dimension of an almost even-Clifford Hermitian manifold, $S$
denotes the homomorphic image of the aforementioned Spin group 
determined by the even-Clifford structure, 
and $\mathcal{N}_{SO(N)}(S)$ denotes 
its normalizer in $SO(N)$.
In fact, we will verify that there is a lift for the connected component of the identity 
$\mathcal{N}^0_{SO(N)}(S)$, 
since the other components are diffeomorphic to it and will also lift to the Spin group. 
Furthermore, note that an almost even-Clifford Hermitian manifold might still be Spin 
even if there is no such a lifting map, as in the case of quaternionic projective spaces $\mathbb{H}\mathbb{P}^m$
of odd quaternionic dimension $m$. It would be interesting, at least for the authors, to find and study 
non-Spin almost even-Clifford manifolds of rank 4, 6 and 8.

The note is organized as follows.
In Section~\ref{sec: preliminaries}, we recall some
preliminaries on Clifford algebras, the Spin group and representations, almost even-Clifford manifolds, etc.
In Section \ref{sec: complexifications}, we determine the complexifications of 
real representations of even 
Clifford algebras containing no trivial summands (cf. Theorem \ref{theo: complexifications}), 
identify the subgroups $\mathcal{N}^0_{SO(N)}(S)$ 
as finite quotients of products of classical groups (or real lines in some cases) and spin groups  (cf. Theorem \ref{theo: quotient groups}),
and calculate their fundamental groups giving explicit generators (cf. Theorem \ref{theo: fundamental groups}).
In Section~\ref{sec: lifts}, we determine when the aformentioned lifts exist (cf. Theorem \ref{theo: lifts}).

\vspace{.1in}

\noi {\em Acknowledgements}.  The third author wishes to thank
the International Centre for Theoretical Physics and 
the Institut des Hautes \'Etudes Scientifiques for their hospitality and support.

\section{Preliminaries
}\label{sec: preliminaries}

The material presented in this section can be consulted in \cite{friedrich}.

\subsection{Clifford algebra, spin group and representation}
Let $Cl_n$ denote the $2^n$-dimensional real Clifford algebra generated by the orthonormal vectors
$e_1, e_2, \ldots, e_n\in \mathbb{R}^n$ subject to
the relations
\[e_i e_j + e_j e_i = -2\delta_{ij},\]
and
$\mathbb{C}l_n=Cl_n\otimes_{\mathbb{R}}\mathbb{C}$ its complexification. 
The even Clifford subalgebra $Cl_r^0$ is defined as the invariant (+1)-subspace of the involution of $Cl_r$
induced 
by the map $-{\rm Id}_{\mathbb{R}^r}$.
For any vector $Y=y_1e_1+\cdots+y_ne_n$, the product
\[e_iYe_i= y_1e_1+\cdots+y_{i-1}e_{i-1}-y_ie_i+y_{i+1}e_{i+1}+\cdots+y_ne_n\]
gives the reflection of the $i$-th coordinate, and the conjugation with the volume element ${\rm vol}_n=e_1\cdots e_n$ gives the reflection
on the origin of $\mathbb{R}^n$, i.e.
\[(e_1\cdots e_n)Y(e_n\cdots e_1) = -Y.\]

There exist algebra isomorphisms
\[\mathbb{C}l_n\cong \left\{
                     \begin{array}{ll}
                     \End(\mathbb{C}^{2^k}) & \mbox{if $n=2k$,}\\
                     \End(\mathbb{C}^{2^k})\oplus\End(\mathbb{C}^{2^k}) & \mbox{if $n=2k+1$,}
                     \end{array}
\right.
\]
and the space of (complex) spinors is defined to be
\[\Delta_n:=\mathbb{C}^{2^k}=\underbrace{\mathbb{C}^2\otimes \ldots \otimes \mathbb{C}^2}_{k\,\,\,\,times}.\]
The map
\[\kappa:\mathbb{C}l_n \lra \End(\mathbb{C}^{2^k})\]
is defined to be either the aforementioned isomorphism for $n$ even, or the isomorphism followed
by the projection onto the first summand for $n$ odd.
In order to make $\kappa$ explicit, consider the following matrices
\[Id = \left(\begin{array}{ll}
1 & 0\\
0 & 1
      \end{array}\right),\quad
g_1 = \left(\begin{array}{ll}
i & 0\\
0 & -i
      \end{array}\right),\quad
g_2 = \left(\begin{array}{ll}
0 & i\\
i & 0
      \end{array}\right),\quad
T = \left(\begin{array}{ll}
0 & -i\\
i & 0
      \end{array}\right).
\]
In terms of the generators $e_1, \ldots, e_n$ of the Clifford algebra, $\kappa$ can be
described explicitly as follows,
\begin{eqnarray}
e_1&\mapsto& Id\otimes Id\otimes \ldots\otimes Id\otimes Id\otimes g_1,\nonumber\\
e_2&\mapsto& Id\otimes Id\otimes \ldots\otimes Id\otimes Id\otimes g_2,\nonumber\\
e_3&\mapsto& Id\otimes Id\otimes \ldots\otimes Id\otimes g_1\otimes T,\nonumber\\
e_4&\mapsto& Id\otimes Id\otimes \ldots\otimes Id\otimes g_2\otimes T,\nonumber\\
\vdots && \dots\nonumber\\
e_{2k-1}&\mapsto& g_1\otimes T\otimes \ldots\otimes T\otimes T\otimes T,\nonumber\\
e_{2k}&\mapsto& g_2\otimes T\otimes\ldots\otimes T\otimes T\otimes T,\nonumber
\end{eqnarray}
and, if $n=2k+1$, 
\[ e_{2k+1}\mapsto i\,\, T\otimes T\otimes\ldots\otimes T\otimes T\otimes T.\]
The vectors 
\[u_{+1}={1\over \sqrt{2}}(1,-i)\quad\quad\mbox{and}\quad\quad u_{-1}={1\over \sqrt{2}}(1,i),\]
form a unitary basis of $\mathbb{C}^2$ with respect to the standard Hermitian product.
Thus, 
\[\mathcal{B}=\{u_{\varepsilon_1,\ldots,\varepsilon_k}=u_{\varepsilon_1}\otimes\ldots\otimes
u_{\varepsilon_k}\,\,|\,\, \varepsilon_j=\pm 1,
j=1,\ldots,k\},\]
is a unitary basis of $\Delta_n=\mathbb{C}^{2^k}$
with respect to the naturally induced Hermitian product.
We will denote inner and Hermitian products (as well as Riemannian and Hermitian
metrics) by the same symbol $\left<\cdot,\cdot\right>$ trusting that the context will make clear
which product is being used.

A quaternionic structure $\alpha$ on $\mathbb{C}^2$ is given by
\[\alpha\left(\begin{array}{c}
z_1\\
z_2
              \end{array}
\right) = \left(\begin{array}{c}
-\overline{z}_2\\
\overline{z}_1
              \end{array}\right),\]
and a real structure $\beta$ on $\mathbb{C}^2$ is given by
\[\beta\left(\begin{array}{c}
z_1\\
z_2
              \end{array}
\right) = \left(\begin{array}{c}
\overline{z}_1\\
\overline{z}_2
              \end{array}\right).\]
Following \cite[p. 31]{friedrich}, the real and quaternionic structures $\gamma_n$  on
$\Delta_n=(\mathbb{C}^2)^{\otimes
[n/2]}$ are built as follows
\[
\begin{array}{cclll}
 \gamma_n &=& (\alpha\otimes\beta)^{\otimes 2k} &\mbox{if $n=8k,8k+1$}& \mbox{(real),} \\
 \gamma_n &=& \alpha\otimes(\beta\otimes\alpha)^{\otimes 2k} &\mbox{if $n=8k+2,8k+3$}&
\mbox{(quaternionic),} \\
 \gamma_n &=& (\alpha\otimes\beta)^{\otimes 2k+1} &\mbox{if $n=8k+4,8k+5$}&\mbox{(quaternionic),} \\
 \gamma_n &=& \alpha\otimes(\beta\otimes\alpha)^{\otimes 2k+1} &\mbox{if $n=8k+6,8k+7$}&\mbox{(real).}
\end{array}
\]

The Spin group $Spin(n)\subset Cl_n$ is the subset 
\[Spin(n) =\{x_1x_2\cdots x_{2l-1}x_{2l}\,\,|\,\,x_j\in\mathbb{R}^n, \,\,
|x_j|=1,\,\,l\in\mathbb{N}\},\]
endowed with the product of the Clifford algebra.
It is a Lie group and its Lie algebra is
\[\mathfrak{spin}(n)=\mbox{span}\{e_ie_j\,\,|\,\,1\leq i< j \leq n\}.\]
The restriction of $\kappa$ to $Spin(n)$ defines the Lie group representation
\[\kappa_n:=\kappa|_{Spin(n)}:Spin(n)\lra GL(\Delta_n),\]
which is, in fact, special unitary. 
We have the corresponding Lie algebra representation
\[\kappa_{n_*}:\mathfrak{spin}(n)\lra \mathfrak{gl}(\Delta_n).\] 
Recall that the Spin group $Spin(n)$ is the universal double cover of $SO(n)$, $n\ge 3$. For $n=2$
we consider $Spin(2)$ to be the connected double cover of $SO(2)$.
The covering map will be denoted by 
\[\lambda_n:Spin(n)\rightarrow SO(n)\subset GL(\mathbb{R}^n).\]
Its differential is given
by $\lambda_{n_*}(e_ie_j) = 2E_{ij}$, where $E_{ij}=e_i^*\otimes e_j - e_j^*\otimes e_i$ is the
standard basis of the skew-symmetric matrices, and $e^*$ denotes the metric dual of the vector $e$.
Furthermore, we will abuse the notation and also denote by $\lambda_n$ the induced representation on
the exterior algebra $\ext^*\mathbb{R}^n$. 

By means of $\kappa$, we have the Clifford multiplication
\begin{eqnarray*}
\mu_n:\mathbb{R}^n\otimes \Delta_n &\lra&\Delta_n\\ 
x \otimes \phi &\mapsto& \mu_n(x\otimes \phi)=x\cdot\phi :=\kappa(x)(\phi) .
\end{eqnarray*}
The Clifford multiplication $\mu_n$ 
 is skew-symmetric with respect to the Hermitian product
\begin{equation*}
\left<x\cdot\phi_1 , \phi_2\right> =\left<\mu_n(x\otimes \phi_1) , \phi_2\right> 
=-\left<\phi_1 , \mu_n(x\otimes \phi_2)\right>
=-\left<\phi_1 , x\cdot \phi_2\right>, 
\end{equation*} 
 is  $Spin(n)$-equivariant and can be extended to a $Spin(n)$-equivariant map 
\begin{eqnarray*}
\mu_n:\ext^*(\mathbb{R}^n)\otimes \Delta_n &\lra&\Delta_n\\ 
\omega \otimes \psi &\mapsto& \omega\cdot\psi.
\end{eqnarray*}

When $n$ is even, we define the following involution
\begin{eqnarray*}
\Delta_n&\longrightarrow& \Delta_n \\ 
\psi &\mapsto& (-i)^{n\over 2}{\rm vol}_n\cdot \psi.
\end{eqnarray*}
The $\pm 1$ eigenspace of this involution is denoted $\Delta_n^\pm$. 
These spaces have equal dimension and are irreducible representations of $Spin(n)$.
Note that our definition differs from the one given in \cite{friedrich} by a $(-1)^{n\over 2}$. 
The reason for this difference
is that we want the spinor $u_{1,\ldots,1}$ to be always positive.
In this case, we will denote the two representations by
\begin{eqnarray*}
\kappa_n^\pm:Spin(n)&\lra& GL(\Delta_n^\pm).
\end{eqnarray*}
Note that while these representations are irreducible, they are not faithful, with kernels isomorphic to 
$\mathbb{Z}_2$ if $n\not =4$.

Now, we summarize some results about real representations of $Cl_r^0$ in the next table (cf. \cite{Lawson}).
Here $d_r$ denotes the dimension of an irreducible representation of $Cl^0_r$ and $v_r$ the number of distinct non-trivial
irreducible representations.
Let $\tilde\Delta_r$ denote the irreducible representation of $Cl_r^0$ for $r\not\equiv0$ $(\mbox{mod } 4) $
and $\tilde\Delta^{\pm}_r$ denote the irreducible representations for $r\equiv0$ $(\mbox{mod } 4)$.
\[\begin{array}{|c|c|c|c|c|c|}
\hline
r \mbox{ (mod 8)}&d_r&Cl_r^0&\tilde\Delta_r\,\,/\,\,\tilde\Delta_r^\pm\cong \mathbb{R}^{d_r}& v_r \tstrut\\
\hline
1&2^{\lfloor\frac r2\rfloor}&\mathbb R(d_r)&\mathbb{R}^{d_r}&1 \tstrut\\
\hline
2&2^{\frac r2}&\mathbb C(d_r/2)&\mathbb{C}^{d_r/2}&1 \tstrut\\
\hline
3&2^{\lfloor\frac r2\rfloor+1}&\mathbb H(d_r/4)&\mathbb{H}^{d_r/4}&1 \tstrut\\
\hline
4&2^{\frac r2}&\mathbb H(d_r/4)\oplus \mathbb H(d_r/4)&\mathbb{H}^{d_r/4}&2 \tstrut\\
\hline
5&2^{\lfloor\frac r2\rfloor+1}&\mathbb H(d_r/4)&\mathbb{H}^{d_r/4}&1 \tstrut\\
\hline
6&2^{\frac r2}&\mathbb C(d_r/2)&\mathbb{C}^{d_r/2}&1 \tstrut\\
\hline
7&2^{\lfloor\frac r2\rfloor}&\mathbb R(d_r)&\mathbb{R}^{d_r}&1 \tstrut\\
\hline
8&2^{\frac r2-1}&\mathbb R(d_r)\oplus \mathbb R(d_r)&\mathbb{R}^{d_r}&2 \tstrut\\
\hline
\end{array}
\]
\centerline{Table 1}

\subsection{Maximal Torus of $Spin(r)$}
In this subsection, we recall explicit expressions for elements of the maximal torus of the Spin group since it will be useful 
to consider paths within such torus.

The rotation
\[\left(
\begin{array}{ccccc}
\cos(\theta_1) & -\sin(\theta_1) &  &  & \\
\sin(\theta_1) & \cos(\theta_1) &   &  & \\
 &  & 1 &    & \\
 &   &  & \ddots & \\
 &    &  &  & 1
\end{array}
\right)_{r\times r}
\]
can be achieved by using the element
\[e_1(-\cos(\theta_1/2)e_1 + \sin(\theta_1/2)e_2)=\cos(\theta_1/2) + \sin(\theta_1/2)e_1e_2 \in Spin(r)\]
as follows 
\[(\cos(\theta_1/2) + \sin(\theta_1/2)e_1e_2)y(\cos(\theta_1/2) - \sin(\theta_1/2)e_1e_2) \]
\[= (y_1\cos(\theta_1)-y_2\sin(\theta_1))e_1
+ (y_1\sin(\theta_1)+y_2\cos(\theta_1)e_2+y_3e_3+\cdots+y_re_r,\]
for $y=y_1e_1+\cdots+y_re_r\in\mathbb{R}^r$. 
Thus, we see that the corresponding elements in  $Spin(r)$ are exactly
\[\pm(\cos(\theta_1/2) + \sin(\theta_1/2)e_1e_2).\]
Furthermore, we can see that a maximal torus of $Spin(r)$ consists of elements of the form
\[t(\theta_1,\ldots,\theta_{[{r\over 2}]})=\prod_{j=1}^{[{r\over 2}]} (\cos(\theta_j/2) + \sin(\theta_j/2)e_{2j-1}e_{2j}),\]
noting that the parameters $\theta_j$ must now run between $0$ and $4\pi$. 
Furthermore, using the explicit description of the isomorphisms given above, we see that
\begin{eqnarray*}
(\cos(\theta_1/2) + \sin(\theta_1/2)e_1e_2)\cdot u_{\varepsilon_1,\ldots\varepsilon_k} 
   &=&  
  \cos(\theta_1/2) u_{\varepsilon_1,\ldots\varepsilon_k} + \sin(\theta_1/2)e_1e_2 \cdot u_{\varepsilon_1,\ldots\varepsilon_k}\\
   &=&  
  \cos(\theta_1/2) u_{\varepsilon_1,\ldots\varepsilon_k} + i\varepsilon_k\sin(\theta_1/2)u_{\varepsilon_1,\ldots\varepsilon_k}\\
   &=&  
  (\cos(\theta_1/2)  + i\varepsilon_k\sin(\theta_1/2))u_{\varepsilon_1,\ldots\varepsilon_k}\\
   &=&  
  e^{i{\varepsilon_k  \theta_1 \over 2}}u_{\varepsilon_1,\ldots\varepsilon_k}, 
\end{eqnarray*}
and similarly,
\[t(\theta_1,\ldots,\theta_{[{r\over 2}]})\cdot u_{\varepsilon_1,\ldots\varepsilon_{[{r\over 2}]}}
=e^{{i\over 2}{\sum_{j=1}^{[{r\over 2}]}\varepsilon_{k+1-j}\theta_j}}\cdot u_{\varepsilon_1,\ldots\varepsilon_k}.\]
Thus, the basis vectors $u_{\varepsilon_1,\ldots\varepsilon_k}$ are weight vectors of the standard spin 
representation with weight 
\[{1\over 2}\sum_{j=1}^{[{r\over 2}]}\varepsilon_{[{r\over 2}]+1-j}\theta_j,\]
which in coordinate vectors are the well known expressions
\[\left(\pm{1\over 2},\pm{1\over 2},\cdots,\pm{1\over 2}\right).\]
Moreover, in terms of the (appropriately ordered) basis $\mathcal{B}$, the matrix associated to an element $t(\theta_1,\ldots,\theta_{[{r\over 2}]})$ is
{\tiny
\[
\left(\begin{array}{ccccccc}
e^{{i\over 2}(\theta_1+\theta_2+\cdots+\theta_{[{r\over 2}]})} &  &  &  &  & &\\
 & e^{{i\over 2}(-\theta_1+\theta_2+\cdots+\theta_{[{r\over 2}]})} &  &  &  & &\\
 &  & e^{{i\over 2}(\theta_1-\theta_2+\cdots+\theta_{[{r\over 2}]})} &  &  & &\\
 &  &  & \ddots & &  & \\
 &  &  & & e^{{i\over 2}(-\theta_1-\theta_2+\cdots+\theta_{[{r\over 2}]})} &  & \\
 & &  &  &  & \ddots & \\
 & &  &  &  &  & e^{{i\over 2}(-\theta_1-\theta_2-\cdots-\theta_{[{r\over 2}]})}
  \end{array}\right).
\]
}
Note that, when $r$ is even, $\Delta_r^+$ is generated by the basis vectors $u_{\varepsilon_1,\ldots\varepsilon_{{r\over 2}}}$ with an even  number  
of $\varepsilon_j$ equal to $-1$, and $\Delta_r^-$ is generated by the basis vectors $u_{\varepsilon_1,\ldots\varepsilon_{{r\over 2}}}$ with an odd number  
of $\varepsilon_j$ equal to $-1$. Therefore, after reordering the basis, the matrix above can be split into two blocks of equal size: 
one block in which the exponents contain an even number of negative signs 
{\tiny
\begin{equation*}
 \left(\begin{array}{ccccc}
e^{{i\over 2}(\theta_1+\theta_2+\cdots+\theta_{{r\over 2}})} &  &  &  &   \\
 & e^{{i\over 2}(-\theta_1-\theta_2+\cdots+\theta_{{r\over 2}})} &  &  &   \\
 &  & e^{{i\over 2}(-\theta_1+\theta_2-\cdots+\theta_{{r\over 2}})} &  &   \\
 &  &  &  & \ddots 
  \end{array}\right)
\end{equation*}
}
and
another block in which the exponents contain an odd number of negative signs
{\tiny
\begin{equation*}
\left(\begin{array}{ccccc}
e^{{i\over 2}(-\theta_1+\theta_2+\cdots+\theta_{{r\over 2}})} &  &  &  &   \\
 & e^{{i\over 2}(\theta_1-\theta_2+\cdots+\theta_{{r\over 2}})} &  &  &   \\
 &  & e^{{i\over 2}(\theta_1+\theta_2-\theta_3+\cdots+\theta_{{r\over 2}})} &  &   \\
 &  &  &  & \ddots  
  \end{array}\right).
\end{equation*}
}

\subsection{Even Clifford structures}

\subsubsection{Linear almost even-Clifford Hermitian structures}

\begin{defi}  Let $N\in \mathbb{N}$ and $(e_1,\ldots,e_r)$ an orthonormal frame of $\mathbb{R}^r$.
\begin{itemize}
\item A {\em linear even-Clifford structure of rank $r$} on $\mathbb{R}^N$ is an algebra
representation
\[\Phi:Cl_r^0\longrightarrow \End(\mathbb{R}^N).\]
\item A {\em linear even-Clifford Hermitian structure of rank $r$} on $\mathbb{R}^N$ 
(endowed with a positive definite inner product) is a linear even-Clifford structure of rank $r$ 
such that each bivector $e_ie_j$, $1\leq i< j \leq r$, is mapped to an skew-symmetric 
endomorphism $\Phi(e_ie_j)=J_{ij}$.
\end{itemize}
\end{defi}

{\bf Remarks}.
\begin{itemize}
\item Note that 
\begin{equation}
J_{ij}^2=-{\rm Id}_{\mathbb{R}^N}.\label{eq: almost complex structures} 
\end{equation}
\item Given a linear even-Clifford structure of rank $r$ on $\mathbb{R}^N$, we can average the standard inner product 
$\left<,\right>$ on $\mathbb{R}^N$ as follows
\[(X,Y)=\sum_{k=1}^{[r/2]} \left[\sum_{1\leq i_1<\ldots<i_{2k}<r}
\left<\Phi(e_{i_1\ldots i_{2k}})(X),\Phi(e_{i_1\ldots i_{2k}})(Y)\right>\right], \]
where $(e_1,\ldots,e_r)$ is an orthonormal frame of $\mathbb{R}^r$,
so that the linear even-Clifford structure is Hermitian with respect to the averaged inner product.
\item 
 Given a linear even-Clifford Hermitian structure structure of rank $r$, 
the subalgebra $\mathfrak{spin}(r)$ is mapped injectively into the skew-symmetric
endomorphisms $\End^-(\mathbb{R}^N)$. 
\end{itemize}

\subsubsection{Branching of $\mathbb{R}^N$}\label{subsubsec: branching}
From now on, we will denote by ${\rm Id}_n$ the identity endomorphism of a real/complex $n$-dimensional vector space.

First, let us assume $r\not\equiv 0
\,\,\,({\rm mod}\,\,\,\, 4)$, $r>1$.
In this case, $\mathbb{R}^N$ decomposes into a sum of irreducible representations of $Cl_r^0$.
Since $Cl_r^0$ is simple, its irreducible representations are either trivial or
the standard representation
$\tilde\Delta_r$ of $Cl_r^0$ (cf. \cite{Lawson}). Due to
\rf{eq: almost complex structures}, there are no trivial summands $\mathbb{R}^N$, i.e.
\begin{equation*}
\mathbb{R}^N = \mathbb{R}^m\otimes \tilde\Delta_r 
\end{equation*}
for some $m\in\mathbb{N}$.
Thus, we see that $\mathfrak{spin}(r)$ has an isomorphic image
\[\widehat{\mathfrak{spin}(r)}:= {\rm Id}_{m}  \otimes \kappa_{r^*}(\mathfrak{spin}(r))\subset
\mathfrak{so}(d_rm).\]

Secondly, let us assume $r\equiv 0
\,\,\,({\rm mod}\,\,\,\, 4)$.
Recall that if $\hat\Delta_r$ is the irreducible representation of $Cl_r$, then by restricting this
representation to $Cl^0_r$ it splits as the sum of two inequivalent irreducible representations
\[\hat\Delta_r = \tilde{\Delta}_r^+ \oplus \tilde{\Delta}_r^-.\]
Since $\mathbb{R}^N$ is a representation of $Cl_r^0$ satisfying
\rf{eq: almost complex structures}, there are no trivial summands in $\mathbb{R}^N$
and
\begin{equation*}
\mathbb{R}^N
=  \mathbb{R}^{m_1}\otimes\tilde\Delta_r^+ \oplus  \mathbb{R}^{m_2}\otimes\tilde\Delta_r^-
\end{equation*}
for some $m_1,m_2\in \mathbb{N}$.
By restricting this representation to $\mathfrak{spin}(r)\subset Cl_r^0$,
consider the isomorphic image
\[
\widehat{\mathfrak{spin}(r)}:=\{  {\rm Id}_{m_1}\otimes \xi^+\oplus
{\rm Id}_{m_2}\otimes \xi^-\, | \,\,\xi\in \mathfrak{spin}(r), \xi^+=\kappa_{r^*}^+(\xi), \xi^-=\kappa_{r^*}^-(\xi)\}.
\]

\subsubsection{Almost even-Clifford Hermitian manifolds}

\begin{defi} Let  $r\geq 2$.
\begin{itemize}  
 \item 
A {\em rank $r$ almost even-Clifford structure} on a smooth manifold $M$
is a smoothly varying choice of a rank $r$ linear even-Clifford structure on each tangent space of $M$.

\item A smooth manifold carrying an almost even-Clifford structure will be called an {\em almost even-Clifford 
manifold}. 

 \item 
A {\em rank $r$ almost even-Clifford Hermitian structure} on a Riemannian manifold
$M$
is a smoothly varying choice of a linear even-Clifford Hermitian structure on each tangent space of $M$.

\item A Riemannian manifold carrying an almost even-Clifford Hermitian structure will be called an {\em almost even-Clifford Hermitian
manifold}.

\end{itemize}
\end{defi}

{\bf Remark}.
Our definition of almost even-Clifford hermintian structure implies that in \cite{Moroianu-Semmelmann}.

\section{Complexifications, structure groups and fundamental groups} \label{sec: complexifications}

In this section we will study 
the complexification $\mathbb{R}^N\otimes \mathbb{C}$ and its decomposition as a representation of 
the normalizers $\mathcal{N}^0_{\mathfrak{so}(N)}(S)$ (cf. \ref{theo: complexifications}), where $S$ is one of the homomorphic images of $\mathfrak{spin}(r)$ described in Subsection \ref{subsubsec: branching}.
Once we have described such complexifications and decompositions, we will compute the (connected components of the identity of the) structure groups determined by linear even-Clifford structures.
They must be closed Lie subgroups of $SO(N)$ whose Lie algebra is the aforementioned normalizer.
They are actually isomorphic to finite quotients of 
products of classical compact Lie groups (or real lines) with $Spin(r)$ (cf. \ref{theo: quotient groups}). Along the way, we will also compute their fundamental groups to be used in the last section.
All of this will be done in case by case analysis.

Along the way, 
we will introduce notation that will enable us to state
Theorems \ref{theo: complexifications}, \ref{theo: quotient groups} and \ref{theo: fundamental groups}.
The main steps in each case are the following:
\begin{enumerate}
 \item Complexification.
 \begin{itemize}
 \item Identify $\mathbb{R}^N\otimes \mathbb{C}$ as a representation  
 \[G\times Spin(r)\xrightarrow{\rho} SO(N)\subset {\rm Aut}(\mathbb{R}^N\times \mathbb{C}),\] 
 where $G$ denotes a (semi-)simple compact Lie group or $SO(2)$.
 \end{itemize}
 \item Structure group.
 \begin{itemize}
 \item Compute the image of $\rho$, ${\rm Im}(\rho)$.
 \item Compute $\ker(\rho)$ to get
 \[{\rm Im}(\rho)\cong {G\times Spin(r)\over \ker(\rho)}.\]
 \end{itemize}
 \item Fundamental group.
 \begin{itemize}
 \item Identify the universal covering $\widetilde{{\rm Im}(\rho)}\xrightarrow{\tilde\rho} {\rm Im}(\rho)$.
 \item Compute the preimage of $\tilde\rho^{-1}({\rm Id}_N)$ to obtain an explicit description of 
 the fundamental group $\pi_1({\rm Im}(\rho))$.
 \end{itemize}
\end{enumerate}

First, let us recall the following table \cite{Arizmendi-Herrera}, whose entries' precise description will be recalled in each case.
\[  
\begin{array}{|c|c|c|c|}
 \hline
  r \mbox{\ {\rm (mod 8)} }  &N
&\mathcal{C}_{ \mathfrak{so}(N)}(\widehat{ \mathfrak{spin}({r})})
&\mathcal{N}^0_{ \mathfrak{so}(N)}(\widehat{ \mathfrak{spin}({r})})
\rule{0pt}{3ex}\\
\hline
 0 &
d_r(m_1+m_2)& \widehat{\mathfrak{so}(m_1)} \oplus \widehat{\mathfrak{so}(m_2)}  
                    & \widehat{\mathfrak{so}(m_1)}  \oplus  \widehat{\mathfrak{so}(m_2)} \oplus \widehat{\mathfrak{spin}(r)} \tstrut\\
 \hline
 1,7 & d_rm& \widehat{\mathfrak{so}(m)} & \widehat{\mathfrak{so}(m)} \oplus \widehat{\mathfrak{spin}(r)}\tstrut\\
\hline
 2,6 & d_rm&  \widehat{\mathfrak{u}(m)}  &  \widehat{\mathfrak{u}(m)}  \oplus\widehat{\mathfrak{spin}(r)}\tstrut\\
\hline
 3,5 & d_rm&  \widehat{\mathfrak{sp}(m)}  &  \widehat{\mathfrak{sp}(m)}  \oplus\widehat{\mathfrak{spin}(r)}\tstrut\\
\hline
 4 & d_r(m_1+m_2)&  \widehat{\mathfrak{sp}(m_1)} \oplus \widehat{\mathfrak{sp}(m_2)} 
                  &  \widehat{\mathfrak{sp}(m_1)}  \oplus \widehat{\mathfrak{sp}(m_2)} \oplus \widehat{\mathfrak{spin}(r)}\tstrut\\
\hline
  \end{array}
\]
\centerline{Table 2}

\subsection{$r\equiv 1,7 \,\,({\rm mod}\,\, 8)$}

\underline{Complexification}.
In this case, $\tilde\Delta_r$ is the subspace of $\Delta_r$ fixed by the corresponding real structure $\gamma_r$, i.e.
\[\tilde\Delta_r \otimes \mathbb{C} = \Delta_r.\]
The centralizer subalgebra of $\widehat{\mathfrak{spin}(r)}$ in $\mathfrak{so}(N)$ is
\[C_{\mathfrak{so}(N)}(\widehat{\mathfrak{spin}(r)})=\mathfrak{so}(m)\otimes {\rm Id}_{\tilde\Delta_r} =: \widehat{\mathfrak{so}(m)},\]
where $N=d_rm$.
If $z\in \mathbb{C}, v\otimes \psi\in\mathbb{R}^m\otimes \tilde\Delta_r$ and $A\in\mathfrak{so}(m)$,
\begin{eqnarray*}
 (A\otimes  {\rm Id}_{\tilde\Delta_r}) (zv\otimes \psi)
   &=& zAv\otimes \psi,
\end{eqnarray*}
which means 
\[(\mathbb{R}^m\otimes \tilde\Delta_r)\otimes \mathbb{C} = \mathbb{C}^m \otimes \Delta_r,\]
where $\mathbb{C}^m$ denotes the standard complex representation of $SO(m)$.
Thus, we have a representation 
\[SO(m)\times Spin(r)\xrightarrow{\rho} SO(N)\subset {\rm Aut}(\mathbb{C}^m \otimes \Delta_r). \]

\noi\underline{Structure group}. 
Since $\widehat{\mathfrak{so}(m)}$ and $\widehat{\mathfrak{spin}(r)}$ commute with each other,
we can take separately the exponentials of their elements within $\mathbb{C}(N)$.
The exponential of 
$A\otimes {\rm Id}_{\tilde\Delta_r}\in \widehat{\mathfrak{so}(m)}$ gives 
\[ e^{A}\otimes{\rm Id}_{\tilde\Delta_r} \in
\widehat{SO(m)}:= SO(m)\otimes{\rm Id}_{\tilde\Delta_r}\cong SO(m).
 \]
On the other hand, if ${\rm Id}_{m}\otimes \xi\in\widehat{\mathfrak{spin}(r)}$, 
its exponential  is
\[{\rm Id}_{m}\otimes e^\xi \in {\rm Id}_{m}\otimes \kappa(Spin(r)) =: \widehat{Spin(r)} \cong Spin(r),\]
since $Spin(r)$ is represented faithfully on $\Delta_r$.
The image of $SO(m)\times Spin(r)$ in $SO(N)$
under the aforementioned representation is 
\[\mathcal{N}_{SO(N)}^0(\widehat{Spin(r)})=\widehat{SO(m)}\widehat{Spin(r)},\] 
the subgroup of all possible products
of elements of the two subgroups, 
i.e. we have 
\[SO(m)\times Spin(r) \xrightarrow{\rho}  \widehat{SO(m)}\widehat{Spin(r)}\subset SO(N).\]
Now we need to find $\ker(\rho)$ and identify $\widehat{SO(m)}\widehat{Spin(r)}$ as a quotient
\[\widehat{SO(m)}\widehat{Spin(r)}\cong {SO(m)\times Spin(r)\over \ker(\rho)} .  \]
If there are elements 
$g\in SO(m)$ and $h \in Spin(r)$ such that
\[\rho(g,h)={\rm Id}_{N},\]
then 
\[\widehat{Spin(r)}\ni\rho({\rm Id}_{m},h)=\rho(g,1)^{-1}\in \widehat{SO(m)}.\]
Since $\rho({\rm Id}_{m},h)$ commutes with every element of $\widehat{Spin(r)}$, it belongs to its center 
$Z(\widehat{Spin(r)})\cong Z(Spin(r)) = \mathbb{Z}_2 =\{\pm 1\}$. Note that  
$-1\in Spin(r)$ maps to $-{\rm Id}_{\Delta_r}$ under the $Spin(r)$ representation $\Delta_r$, and $({\rm Id}_{m},-1)$ 
maps to $-{\rm Id}_{m}\otimes{\rm Id}_{\Delta_r}\in SO(N) $ under $\rho$.
Moreover,  $-{\rm Id}_{m}\otimes{\rm Id}_{\Delta_r}$ belongs to $\widehat{SO(m)}$ only if $m$ is even.
Thus, $\ker(\rho)= \{\pm( {\rm Id}_{m},1 )\}=\mathbb{Z}_2$ and
\[\widehat{SO(m)}\widehat{Spin(r)} \cong {SO(m)\times Spin(r)\over \mathbb{Z}_2}  \]
if $m$ is even, 
and $\ker(\rho)= \{( {\rm Id}_{m},1 )\}$,
\[\widehat{SO(m)}\widehat{Spin(r)} \cong SO(m)\times Spin(r)  \]
if $m$ is odd.

\noi \underline{Fundamental group}. Clearly, we only need to deal with the case when $m$ is even. 
\begin{itemize}
 \item If $m\geq 4$,
let $\tilde\rho$ denote the following composition 
\[\begin{array}{c}
Spin(m)\times Spin(r)\\
\downarrow\\
SO(m)\times Spin(r)\\
\downarrow\\
SO(m)\times_{\mathbb{Z}_2}Spin(r)
  \end{array}
\]
We need to find all the elements of $Spin(m)\times Spin(r)$ that map to
\[\pm( {\rm Id}_{m},1 ).\]
The elements of $Spin(m)\times Spin(r)$ that map to
$( {\rm Id}_{m},1 )\in SO(m)\times Spin(r)$ are
\[ (\pm 1,1),\]
and the elements of $Spin(m)\times Spin(r)$ that map to
$(- {\rm Id}_{m},-1 )\in SO(m)\times Spin(r)$ are
\[(\pm {\rm vol}_m,-1), \] 
i.e.
\begin{eqnarray*}
\ker(\tilde\rho)&=&\{( 1,1),(-1,1),({\rm vol}_m,-1),(- {\rm vol}_m,-1)\}.\\
\pi_1(\widehat{SO(m)}\widehat{Spin(r)})
&\cong&\left\{\begin{array}{ll}
 \mathbb{Z}_2\oplus\mathbb{Z}_2 & \mbox{if $m\equiv 0 \,\,(\mbox{mod 4})$,} \\
 \mathbb{Z}_4 & \mbox{if $m\equiv 2 \,\,(\mbox{mod 4})$.}
                                                 \end{array}
\right. 
\end{eqnarray*}

\item If $m=2$,
let $\tilde\rho$ denote the following composition 
\[\begin{array}{c}
\mathbb{R}\times Spin(r)\\
\downarrow\\
SO(2)\times Spin(r)\\
\downarrow\\
SO(2)\times_{\mathbb{Z}_2}Spin(r)
  \end{array}
\]
Similarly,
\begin{eqnarray*}
\ker(\tilde\rho)&=&\{(2k\pi ,1) \,|\, k\in \mathbb{Z}\}\cup\{((2k+1)\pi,-1)\,|\, k\in \mathbb{Z}\},\\
\pi_1(\widehat{SO(2)}\widehat{Spin(r)})&\cong &\mathbb{Z}. 
\end{eqnarray*}

\end{itemize}

\subsection{$r\equiv 0 \,\,({\rm mod}\,\, 8)$}

\underline{Complexification}. 
In this case, $\tilde\Delta_r^+$ and $\tilde\Delta_r^-$ are the subspaces of $\Delta_r^+$ and $\Delta_r^-$ fixed by the corresponding real structure $\gamma_r$ mentioned in 
Section \ref{sec: preliminaries}, and
\begin{eqnarray*}
\tilde\Delta_r^+ \otimes \mathbb{C} &=& \Delta_r^+,\\
\tilde\Delta_r^- \otimes \mathbb{C} &=& \Delta_r^-. 
\end{eqnarray*}
The centralizer subalgebra of $\widehat{\mathfrak{spin}(r)}$ is 
\begin{eqnarray*}
C_{\mathfrak{so}(N)}(\widehat{\mathfrak{spin}(r)})
&=&\widehat{\mathfrak{so}(m_1)}\oplus\widehat{\mathfrak{so}(m_2)}\\
&=&
\left(
\begin{array}{cc}
\widehat{\mathfrak{so}(m_1)} & \\
 & \widehat{\mathfrak{so}(m_2)} 
\end{array}
\right)\\
&:=&
\left(
\begin{array}{cc}
 \mathfrak{so}(m_1)\otimes {\rm Id}_{\tilde\Delta_r^+} & \\
 & \mathfrak{so}(m_2)\otimes  {\rm Id}_{\tilde\Delta_r^-} 
\end{array}
\right), 
\end{eqnarray*}
where $N=d_r(m_1+m_2)$.
If $A_1\in\mathfrak{so}(m_1)$, $A_2\in\mathfrak{so}(m_2)$, $z_1,z_2\in\mathbb{C}$,  
$v_1\otimes \psi_1  + v_2\otimes\psi_2\in\mathbb{R}^{m_1}\otimes \tilde\Delta_r^+\oplus \mathbb{R}^{m_1}\otimes \tilde\Delta_r^-$,
\[
\left(
\begin{array}{cc}
 A_1\otimes {\rm Id}_{\tilde\Delta_r^+} & \\
 & A_2\otimes  {\rm Id}_{\tilde\Delta_r^-} 
\end{array}
\right)
\left(
\begin{array}{l}
z_1v_1\otimes \psi_1\\
z_2v_2\otimes \psi_2
\end{array}
\right)
=\left(
\begin{array}{l}
z_1A_1v_1\otimes \psi_1\\
z_2A_2v_2\otimes \psi_2
\end{array}
\right),
\]
which means  
\[(\mathbb{R}^{m_1}\otimes \tilde\Delta_r^+\oplus\mathbb{R}^{m_2}\otimes \tilde\Delta_r^-)\otimes \mathbb{C} 
= \mathbb{C}^{m_1} \otimes \Delta_r^+\oplus\mathbb{C}^{m_2} \otimes \Delta_r^-,\]
where $\mathbb{C}^{m_1}$ and $\mathbb{C}^{m_2}$ denote the standard complex representation of $\mathfrak{so}(m_1)$ and $\mathfrak{so}(m_2)$ respectively.
Thus, we have a representation 
\[SO(m_1)\times SO(m_2)\times Spin(r)\longrightarrow SO(N)\subset 
 {\rm Aut}(\mathbb{C}^{m_1}\otimes \Delta_r^+ \oplus \mathbb{C}^{m_2}\otimes \Delta_r^-)
.\]

\noi\underline{Structure group}.
Since $\widehat{\mathfrak{so}(m_1)}\oplus \widehat{\mathfrak{so}(m_2)}$ and $\widehat{\mathfrak{spin}(r)}$ commute, 
we can take the exponentials of their elements separately within $\mathbb{C}(N)$.
The exponential of 
\[\left(
\begin{array}{cc}
 A_1\otimes {\rm Id}_{\Delta_r^+} & \\
 & A_2\otimes  {\rm Id}_{\Delta_r^-} 
\end{array}
\right)\in \widehat{\mathfrak{so}(m_1)}\oplus\widehat{\mathfrak{so}(m_2)}
\]
is
\begin{eqnarray*}
\left(
\begin{array}{cc}
 e^{A_1}\otimes {\rm Id}_{\Delta_r^+} & \\
 & e^{A_2}\otimes  {\rm Id}_{\Delta_r^-} 
\end{array}
\right)
&\in& 
\left(\begin{array}{ll}
\widehat{SO(m_1)} &   \\
 &  \widehat{SO(m_2)} 
\end{array}\right)
\\
&:=&
\left(\begin{array}{ll}
SO(m_1)\otimes{\rm Id}_{\Delta_r^+} &   \\
 &  SO(m_2)\otimes{\rm Id}_{\Delta_r^-} 
\end{array}\right).
\\
&=&
 \widehat{SO(m_1)}\times \widehat{SO(m_2)}.
\end{eqnarray*}
On the other hand, if ${\rm Id}_{m_1}\otimes \xi^+\oplus {\rm Id}_{m_2}\otimes \xi^-\in\widehat{\mathfrak{spin}(r)}$, 
its exponential  is
\begin{eqnarray*}
\left(
\begin{array}{cc}
{\rm Id}_{m_1}\otimes e^{\xi^+} & \\
 & {\rm Id}_{m_2}\otimes e^{\xi^-}
\end{array}
\right)
  &\in& 
  \left(
\begin{array}{cc}
{\rm Id}_{m_1}\otimes \kappa^+(Spin(r)) & \\
 & {\rm Id}_{m_2}\otimes \kappa^-(Spin(r))
\end{array}\right)\\
&=& 
\left\{
\begin{array}{ll}
\widehat{Spin(r)}\cong Spin(r) & \mbox{if $m_1>0$ and $m_2>0$,}\\
\widehat{Spin(r)^+}\cong \kappa^+(Spin(r))& \mbox{if $m_1>0$ and $m_2=0$,}\\
\widehat{Spin(r)^-}\cong \kappa^-(Spin(r)) & \mbox{if $m_1=0$ and $m_2>0$,}
\end{array}
\right.
\end{eqnarray*}
where the first representation of $Spin(r)$ is faithful and the last two are not, with 
\begin{eqnarray*}
 \widehat{Spin(r)^\pm}&\cong& \kappa^\pm(Spin(r)) \cong {Spin(r)\over \{1, \pm{\rm vol}_r\} }.
\end{eqnarray*}
The image of $SO(m_1)\times SO(m_2)\times Spin(r)$ in $SO(N)$
under the aforementioned representations are 
\begin{eqnarray*}
\mathcal{N}_{SO(N)}^0(\widehat{Spin(r)}) &=& (\widehat{SO(m_1)}\times \widehat{SO(m_2)})\widehat{Spin(r)},\\
\mathcal{N}_{SO(N)}^0(\widehat{Spin(r)^+}) &=& \widehat{SO(m_1)}\widehat{Spin(r)^+}, \\
\mathcal{N}_{SO(N)}^0(\widehat{Spin(r)^-}) &=&\widehat{SO(m_2)}\widehat{Spin(r)^-},  
\end{eqnarray*}
respectively,
i.e. in each case we have a map
\begin{eqnarray*}
SO(m_1)\times SO(m_2)\times Spin(r) &\xrightarrow{\rho}&  (\widehat{SO(m_1)}\times \widehat{SO(m_2)})\widehat{Spin(r)}\subset SO(N),\\ 
SO(m_1)\times Spin(r) &\xrightarrow{\rho}&  \widehat{SO(m_1)}\widehat{Spin(r)^-}\subset SO(N),\\ 
SO(m_2)\times Spin(r) &\xrightarrow{\rho} & \widehat{SO(m_2)}\widehat{Spin(r)^+}\subset SO(N). 
\end{eqnarray*}
Now we need to find $\ker(\rho)$ in each case to identify the relevant group as a quotient.
\begin{itemize}
 \item Case $m_1, m_2> 0$.
 If there are elements 
$g_i\in SO(m_i)$ and $h \in Spin(r)$ such that
\[\rho(g_1,g_2,h)={\rm Id}_{N},\]
then 
\[\widehat{Spin(r)}\ni\rho({\rm Id}_{m_1},{\rm Id}_{m_2},h)=\rho(g_1,g_2,1)^{-1}\in \widehat{SO(m_1)}\times\widehat{SO(m_2)}.\]
Since 
$\rho({\rm Id}_{m_1},{\rm Id}_{m_2},h)$ commutes with every element of $\widehat{Spin(r)}$, it belongs to its center 
$Z(\widehat{Spin(r)})\cong Z(Spin(r)) = \{1,-1, {\rm vol}_r,-{\rm vol}_r\}
\cong  \mathbb{Z}_2 \oplus \mathbb{Z}_2 $.
Note that 
\begin{itemize}
 \item the element $-1$ is mapped to $-{\rm Id}_{\Delta_r^\pm}$ in the $Spin(r)$ representations $\Delta_r^\pm$, 
 and $({\rm Id}_{m_1},{\rm Id}_{m_2},-1)$ maps to 
$-({\rm Id}_{m_1}\otimes{\rm Id}_{\Delta_r^+}\oplus {\rm Id}_{m_2}\otimes{\rm Id}_{\Delta_r^-})\in SO(N) $; 
the element $-({\rm Id}_{m_1}\otimes{\rm Id}_{\Delta_r^+}\oplus {\rm Id}_{m_2}\otimes{\rm Id}_{\Delta_r^-})$ belongs to 
$\widehat{SO(m_1)}\times \widehat{SO(m_2)}$ if $m_1\equiv m_2\equiv 0\,\,(\mbox{mod $2$})$; 
\item the element ${\rm vol}_r$ is mapped to $\pm {\rm Id}_{\Delta_r^\pm}$, 
 and  $({\rm Id}_{m_1},{\rm Id}_{m_2},{\rm vol}_r)$   maps to 
$({\rm Id}_{m_1}\otimes{\rm Id}_{\Delta_r^+}\oplus (-1){\rm Id}_{m_2}\otimes{\rm Id}_{\Delta_r^-})\in SO(N) $;
the element
$({\rm Id}_{m_1}\otimes{\rm Id}_{\Delta_r^+}\oplus (-1){\rm Id}_{m_2}\otimes{\rm Id}_{\Delta_r^-})$
belongs to 
$\widehat{SO(m_1)}\times\widehat{SO(m_2)}$ if $m_2\equiv 0\,\,(\mbox{mod $2$})$;
\item the element $-{\rm vol}_r$ is mapped to $\mp {\rm Id}_{\Delta_r^\pm}$, 
 and $({\rm Id}_{m_1},{\rm Id}_{m_2},-{\rm vol}_r)$ maps to 
$((-1){\rm Id}_{m_1}\otimes{\rm Id}_{\Delta_r^+}\oplus {\rm Id}_{m_2}\otimes{\rm Id}_{\Delta_r^-})\in SO(N) $;
the element
$((-1){\rm Id}_{m_1}\otimes{\rm Id}_{\Delta_r^+}\oplus {\rm Id}_{m_2}\otimes{\rm Id}_{\Delta_r^-})$
belongs to 
$\widehat{SO(m_1)}\times \widehat{SO(m_2)}$ if $m_1\equiv 0\,\,(\mbox{mod $2$})$.
\end{itemize}
Thus,
\begin{enumerate}
 \item[(1)] if $m_1\equiv m_2\equiv 1\,\,(\mbox{mod $2$})$,
 \begin{eqnarray*}
\ker(\rho)&=&\{({\rm Id}_{m_1},{\rm Id}_{m_2},1)\},\\
(\widehat{SO(m_1)}\times \widehat{SO(m_2)})\widehat{Spin(r)}&\cong& {SO(m_1)\times SO(m_2)\times Spin(r)} ; 
 \end{eqnarray*}
 \item[(2)] if $m_1\equiv 0\,\,(\mbox{mod $2$})$, $ m_2\equiv 1\,\,(\mbox{mod $2$})$,
 \begin{eqnarray*}
\ker(\rho)&=&\{({\rm Id}_{m_1},{\rm Id}_{m_2},1),(-{\rm Id}_{m_1},{\rm Id}_{m_2},-{\rm vol}_r)\}\cong \mathbb{Z}_2,\\
(\widehat{SO(m_1)}\times \widehat{SO(m_2)})\widehat{Spin(r)}&\cong& {SO(m_1)\times SO(m_2)\times Spin(r)\over \mathbb{Z}_2} ; 
 \end{eqnarray*}
 \item[(3)] if $m_1\equiv 1\,\,(\mbox{mod $2$})$, $ m_2\equiv 0\,\,(\mbox{mod $2$})$,
\begin{eqnarray*}
\ker(\rho)&=&\{({\rm Id}_{m_1},{\rm Id}_{m_2},1),({\rm Id}_{m_1},-{\rm Id}_{m_2},{\rm vol})\}\cong \mathbb{Z}_2,\\ 
(\widehat{SO(m_1)}\times \widehat{SO(m_2)})\widehat{Spin(r)}&\cong &{SO(m_1)\times SO(m_2)\times Spin(r)\over \mathbb{Z}_2} ; 
\end{eqnarray*}
 \item[(4)] if $m_1\equiv  m_2\equiv 0\,\,(\mbox{mod $2$})$,
 \begin{eqnarray*}
 \ker(\rho)&=&\{({\rm Id}_{m_1},{\rm Id}_{m_2},1),
     (-{\rm Id}_{m_1}, -{\rm Id}_{m_2},-1),\\
     &&({\rm Id}_{m_1},-{\rm Id}_{m_2},{\rm vol}),
     (-{\rm Id}_{m_1},{\rm Id}_{m_2},-{\rm vol})
     \}\\
     &\cong& \mathbb{Z}_2\oplus \mathbb{Z}_2,\\
(\widehat{SO(m_1)}\times \widehat{SO(m_2)})\widehat{Spin(r)}&\cong& {SO(m_1)\times SO(m_2)\times Spin(r)\over \mathbb{Z}_2\oplus \mathbb{Z}_2} . 
 \end{eqnarray*}
\end{enumerate}

 \item Case $m_1> 0$, $m_2= 0$.
 If there are elements 
$g_1\in SO(m_1)$ and $h \in Spin(r)$ such that
\[\rho(g_1,h)={\rm Id}_{N},\]
then 
\[\widehat{Spin(r)^+}\ni\rho({\rm Id}_{m_1},h)=\rho(g_1,1)^{-1}\in \widehat{SO(m_1)}.\]
Since 
$\rho({\rm Id}_{m_1},h)$ commutes with every element of $\widehat{Spin(r)^+}$, it belongs to its center 
$Z(\widehat{Spin(r)^+})\cong Z(Spin(r)/\{1,{\rm vol}_r\}) = \{1,-1, {\rm vol}_r,-{\rm vol}_r\}/\{1,{\rm vol}_r\}
\cong  \mathbb{Z}_2  $.
Note that 
\begin{itemize}
 \item the element $-1$ is mapped to $-{\rm Id}_{\Delta_r^+}$ in the $Spin(r)$ representation $\Delta_r^+$, 
 and $({\rm Id}_{m_1},-1)$ maps to 
$-{\rm Id}_{m_1}\otimes{\rm Id}_{\Delta_r^+}\in SO(N) $; the element $-{\rm Id}_{m_1}\otimes{\rm Id}_{\Delta_r^+}$ belongs to 
$\widehat{SO(m_1)}$ if $m_1\equiv  0\,\,(\mbox{mod $2$})$.
\end{itemize}
Thus,
\begin{enumerate}
 \item[(5)] if $m_1\equiv  1\,\,(\mbox{mod $2$})$,
 \begin{eqnarray*}
\ker(\rho)&=&\{({\rm Id}_{m_1},1),({\rm Id}_{m_1},{\rm vol}_r)\},\\
\widehat{SO(m_1)}\widehat{Spin(r)^+}&\cong& {SO(m_1)\times Spin(r)\over \mathbb{Z}_2} ; 
 \end{eqnarray*}
 \item[(6)] if $m_1\equiv 0\,\,(\mbox{mod $2$})$,
 \begin{eqnarray*}
\ker(\rho)&=&\{({\rm Id}_{m_1},1),(-{\rm Id}_{m_1},-1),({\rm Id}_{m_1},{\rm vol}_r),(-{\rm Id}_{m_1},-{\rm vol}_r)\},\\
\widehat{SO(m_1)}\widehat{Spin(r)^+}&\cong& {SO(m_1)\times Spin(r)\over \mathbb{Z}_2\oplus \mathbb{Z}_2} .
 \end{eqnarray*}
\end{enumerate}

 \item Case $m_1= 0$, $m_2> 0$.
 If there are elements 
$g_2\in SO(m_2)$ and $h \in Spin(r)$ such that
\[\rho(g_2,h)={\rm Id}_{N},\]
then 
\[\widehat{Spin(r)^-}\ni\rho({\rm Id}_{m_2},h)=\rho(g_2,1)^{-1}\in \widehat{SO(m_2)}.\]
Since 
$\rho({\rm Id}_{m_2},h)$ commutes with every element of $\widehat{Spin(r)^-}$, it belongs to its center 
$Z(\widehat{Spin(r)^-})\cong Z(Spin(r)/\{1,-{\rm vol}_r\}) = \{1,-1, {\rm vol}_r,-{\rm vol}_r\}/\{1,-{\rm vol}_r\}
\cong  \mathbb{Z}_2  $.
Note that 
\begin{itemize}
 \item the element $-1$ is mapped to $-{\rm Id}_{\Delta_r^-}$ in the $Spin(r)$ representation $\Delta_r^-$, 
 and $({\rm Id}_{m_2},-1)$ maps to 
$-{\rm Id}_{m_2}\otimes{\rm Id}_{\Delta_r^-}\in SO(N) $; 
 the element $-{\rm Id}_{m_2}\otimes{\rm Id}_{\Delta_r^-}$ belongs to 
$\widehat{SO(m_2)}$ if $m_2\equiv 0\,\,(\mbox{mod $2$})$;
\end{itemize}
Thus,
\begin{enumerate}
 \item[(7)] if $ m_2\equiv 1\,\,(\mbox{mod $2$})$,
 \begin{eqnarray*}
\ker(\rho)&=&\{({\rm Id}_{m_2},1),({\rm Id}_{m_2},-{\rm vol}_r)\},\\
\widehat{SO(m_2)}\widehat{Spin(r) ^-}&\cong& {SO(m_2)\times Spin(r)\over \mathbb{Z}_2} ; 
 \end{eqnarray*}
 \item[(8)] if $ m_2\equiv 0\,\,(\mbox{mod $2$})$,
 \begin{eqnarray*}
\ker(\rho)&=&\{({\rm Id}_{m_2},1),(-{\rm Id}_{m_2},-1),(-{\rm Id}_{m_2},{\rm vol}_r),({\rm Id}_{m_2},-{\rm vol}_r)\},\\
\widehat{SO(m_2)}\widehat{Spin(r)^-}&\cong& {SO(m_2)\times Spin(r)\over \mathbb{Z}_2\oplus \mathbb{Z}_2} .
 \end{eqnarray*}
\end{enumerate}

\end{itemize}

\noi\underline{Fundamental group}.
We will now analyze each of the previous eight cases:
\begin{itemize}
 \item[(1)] Recall that $m_1,m_2>0$ and $m_1\equiv m_2\equiv 1$ (mod 2).

\begin{itemize}
\item If $m_1,m_2\geq 3$, 
let $\tilde\rho$ denote the following map
 \[\begin{array}{c}
Spin(m_1)\times Spin(m_2)\times Spin(r)\\
\downarrow\\
SO(m_1)\times SO(m_2)\times Spin(r)
   \end{array}
\]
Thus
\begin{eqnarray*}
\ker(\tilde\rho) &=&  \{(1,1,1), (-1,1,1), (1,-1,1),  (-1,-1,1)\},\\    
\pi_1((\widehat{SO(m_1)}\times \widehat{SO(m_2)})\widehat{Spin(r)})&\cong& \mathbb{Z}_2\oplus \mathbb{Z}_2.
\end{eqnarray*}

\item If $m_1=1,m_2\geq 3$, $SO(m_1)=\{{\rm Id}_{1}\}$ and
let $\tilde\rho$ denote the following map
 \[\begin{array}{c}
\{{\rm Id}_{1}\}\times Spin(m_2)\times Spin(r)\\
\downarrow\\
\{{\rm Id}_{1}\}\times SO(m_2)\times Spin(r)
   \end{array}
\]
Thus
\begin{eqnarray*}
\ker(\tilde\rho)&=&\{({\rm Id}_{1},1,1),  ({\rm Id}_{1},-1,1)  \},\\
\pi_1((\widehat{SO(1)}\times \widehat{SO(m_2)})\widehat{Spin(r)})
&\cong& \mathbb{Z}_2.  
\end{eqnarray*}

\item If $m_1\geq 3,m_2=1$,  
let $\tilde\rho$ denote the following map
 \[\begin{array}{c}
Spin(m_1)\times \{{\rm Id}_{1}\}\times Spin(r)\\
\downarrow\\
SO(m_1)\times \{{\rm Id}_{1}\}\times Spin(r)
   \end{array}
\]
Thus
\begin{eqnarray*}
\ker(\tilde\rho)&=&\{(1,{\rm Id}_{1},1), (-1,{\rm Id}_{1},1) \},\\ 
\pi_1((\widehat{SO(m_1)}\times \widehat{SO(1)})\widehat{Spin(r)})&\cong& \mathbb{Z}_2.  
\end{eqnarray*}

\item If $m_1=1,m_2=1$, 
let $\tilde\rho$ be
 \[\begin{array}{c}
\{{\rm Id}_{1}\}\times \{{\rm Id}_{1}\}\times Spin(r)\\
\downarrow\\
\{{\rm Id}_{1}\}\times \{{\rm Id}_{1}\}\times Spin(r)
   \end{array}
\]
Thus,
\begin{eqnarray*}
\ker(\tilde\rho)&=&\{({\rm Id}_{1},{\rm Id}_{1},1)\},\\
\pi_1((\widehat{SO(1)}\times \widehat{SO(1)})\widehat{Spin(r)})&\cong& \{1\}. 
\end{eqnarray*}

\end{itemize} 
 
 \item[(2)] 
Recall that $m_1,m_2>0$, $m_1\equiv 0\,\,(\mbox{mod $2$})$, $ m_2\equiv 1\,\,(\mbox{mod $2$})$.
\begin{itemize}
\item If $m_1\geq 4, m_2\geq 3$,
let $\tilde\rho$ denote the composition
 \[\begin{array}{c}
Spin(m_1)\times Spin(m_2)\times Spin(r)\\
\downarrow\\
SO(m_1)\times SO(m_2)\times Spin(r)\\
\downarrow\\
{SO(m_1)\times SO(m_2)\times Spin(r)\over \mathbb{Z}_2}
   \end{array}
\]
Thus
\begin{eqnarray*}
\ker(\tilde\rho)
   &=&
   \left\{   \begin{array}{ll}
\left<(1,-1,1),(-1,1,1),({\rm vol}_{m_1},1,-{\rm vol}_r)\right> & \mbox{if $m_1\equiv 0$ (mod 4),} \\
\left<(1,-1,1),({\rm vol}_{m_1},1,-{\rm vol}_r)\right> & \mbox{if $m_1\equiv 2$ (mod 4),}
                                              \end{array}
\right.\\ 
\pi_1((\widehat{SO(m_1)}\times \widehat{SO(m_2)})\widehat{Spin(r)})
   &\cong&
   \left\{   \begin{array}{ll}
\mathbb{Z}_2\oplus\mathbb{Z}_2\oplus\mathbb{Z}_2 & \mbox{if $m_1\equiv 0$ (mod 4),} \\
\mathbb{Z}_2\oplus\mathbb{Z}_4 & \mbox{if $m_1\equiv 2$ (mod 4).}
                                              \end{array}
\right.
\end{eqnarray*}

\item If $m_1=2 , m_2\geq 3$, $SO(m_1)=SO(2)$
and let $\tilde\rho$ denote the following composition
 \[\begin{array}{c}
\mathbb{R}\times Spin(m_2)\times Spin(r)\\
\downarrow\\
SO(2)\times SO(m_2)\times Spin(r)\\
\downarrow\\
{SO(2)\times SO(m_2)\times Spin(r)\over \mathbb{Z}_2}
   \end{array}
\]
Thus,
\begin{eqnarray*}
\ker(\tilde\rho) &=& \left<(0,-1,1),(\pi,1,-{\rm vol}_r)\right>,\\
\pi_1((\widehat{SO(2)}\times \widehat{SO(m_2)})\widehat{Spin(r)})
   &=& \mathbb{Z}_2\oplus\mathbb{Z}.
\end{eqnarray*}

\item If $m_1\geq 4, m_2= 1$, 
let $\tilde\rho$ denote the composition
 \[\begin{array}{c}
Spin(m_1)\times \{{\rm Id}_{1}\}\times Spin(r)\\
\downarrow\\
SO(m_1)\times \{{\rm Id}_{1}\}\times Spin(r)\\
\downarrow\\
{SO(m_1)\times \{{\rm Id}_{1}\}\times Spin(r)\over \mathbb{Z}_2}
   \end{array}
\]
Thus
\begin{eqnarray*}
\ker(\tilde\rho)
   &=&
   \left\{   \begin{array}{ll}
 \left<(-1,{\rm Id}_{1},1),({\rm vol}_{m_1},{\rm Id}_{1},-{\rm vol}_r)\right> & \mbox{if $m_1\equiv 0$ (mod 4),} \\
 \left<({\rm vol}_{m_1},{\rm Id}_{1},-{\rm vol}_r)\right> & \mbox{if $m_1\equiv 2$ (mod 4),}
                                              \end{array}
\right.\\ 
\pi_1((\widehat{SO(m_1)}\times\widehat{SO(1)})\widehat{Spin(r)})
   &=&
   \left\{   \begin{array}{ll}
\mathbb{Z}_2\oplus\mathbb{Z}_2 & \mbox{if $m_1\equiv 0$ (mod 4),} \\
\mathbb{Z}_4 & \mbox{if $m_1\equiv 2$ (mod 4).}
                                              \end{array}
\right.\\ 
\end{eqnarray*}

\item If $m_1=2, m_2= 1$, 
let $\tilde\rho$ denote the composition 
 \[\begin{array}{c}
\mathbb{R}\times \{{\rm Id}_{1}\}\times Spin(r)\\
\downarrow\\
SO(2)\times \{{\rm Id}_{1}\}\times Spin(r)\\
\downarrow\\
{SO(2)\times \{{\rm Id}_{1}\}\times Spin(r)\over \mathbb{Z}_2}
   \end{array}
\]
Thus
\begin{eqnarray*}
\ker(\tilde\rho)
   &=& \left<(\pi,{\rm Id}_{1},-{\rm vol}_r)\right>,\\
\pi_1((\widehat{SO(2)}\times\widehat{SO(1)})\widehat{Spin(r)})
   &\cong& \mathbb{Z}. 
\end{eqnarray*}

\end{itemize}

 \item[(3)] 
Recall that $m_1,m_2>0$, $m_1\equiv 1\,\,(\mbox{mod $2$})$, $ m_2\equiv 0\,\,(\mbox{mod $2$})$.
\begin{itemize}
\item If $m_1\geq 3, m_2\geq 4$,
let $\tilde\rho$ denote the composition
 \[\begin{array}{c}
Spin(m_1)\times Spin(m_2)\times Spin(r)\\
\downarrow\\
SO(m_1)\times SO(m_2)\times Spin(r)\\
\downarrow\\
{SO(m_1)\times SO(m_2)\times Spin(r)\over \mathbb{Z}_2}
   \end{array}
\]
Thus,
\begin{eqnarray*}
\ker(\tilde\rho)
   &=&
   \left\{   \begin{array}{ll}
\left<(-1,1,1),(1,-1,1),(1,{\rm vol}_{m_2},{\rm vol}_r)\right> & \mbox{if $m_2\equiv 0$ (mod 4),} \\
\left<(-1,1,1),(1,{\rm vol}_{m_2},{\rm vol}_r)\right> & \mbox{if $m_2\equiv 2$ (mod 4),}
                                              \end{array}
\right.\\ 
\pi_1((\widehat{SO(m_1)}\times \widehat{SO(m_2)})\widehat{Spin(r)})
   &=&
   \left\{   \begin{array}{ll}
\mathbb{Z}_2\oplus\mathbb{Z}_2\oplus\mathbb{Z}_2 & \mbox{if $m_2\equiv 0$ (mod 4),} \\
\mathbb{Z}_2\oplus\mathbb{Z}_4 & \mbox{if $m_2\equiv 2$ (mod 4).}
                                              \end{array}
\right.
\end{eqnarray*}

\item If $m_1\geq 3, m_2=2$,
let $\tilde\rho$ denote the composition
 \[\begin{array}{c}
Spin(m_1)\times \mathbb{R}\times Spin(r)\\
\downarrow\\
SO(m_1)\times SO(2)\times Spin(r)\\
\downarrow\\
{SO(m_1)\times SO(2)\times Spin(r)\over \mathbb{Z}_2}
   \end{array}
\]
Thus,
\begin{eqnarray*}
\ker(\tilde\rho)
   &=&\left<(-1,0,1),(1,\pi,{\rm vol}_r)\right>,
\\
\pi_1((\widehat{SO(m_1)}\times \widehat{SO(2)})\widehat{Spin(r)})
   &\cong& \mathbb{Z}_2\oplus \mathbb{Z}.
\end{eqnarray*}

\item If $m_1= 1, m_2\geq 4$, 
let $\tilde\rho$ denote the composition
 \[\begin{array}{c}
\{{\rm Id}_{1}\}\times Spin(m_2)\times Spin(r)\\
\downarrow\\
\{{\rm Id}_{1}\}\times SO(m_2)\times Spin(r)\\
\downarrow\\
{\{{\rm Id}_{1}\}\times SO(m_2)\times Spin(r)\over \mathbb{Z}_2}
   \end{array}
\]
Thus,
\begin{eqnarray*}
\ker(\tilde\rho)
   &=&
   \left\{   \begin{array}{ll}
\left<({\rm Id}_{1},-1,1),({\rm Id}_{1},{\rm vol}_{m_2},{\rm vol}_r)\right> & \mbox{if $m_2\equiv 0$ (mod 4),} \\
\left<({\rm Id}_{1},{\rm vol}_{m_2},{\rm vol}_r)\right> & \mbox{if $m_2\equiv 2$ (mod 4),}
                                              \end{array}
\right.\\ 
\pi_1((\widehat{SO(m_1)}\times\widehat{SO(m_2)})\widehat{Spin(r)})
   &=&
   \left\{   \begin{array}{ll}
\mathbb{Z}_2\oplus\mathbb{Z}_2 & \mbox{if $m_2\equiv 0$ (mod 4),} \\
\mathbb{Z}_4 & \mbox{if $m_2\equiv 2$ (mod 4).}
                                              \end{array}
\right.
\end{eqnarray*}

\item If $m_1=1, m_2=2 $,  
let $\tilde\rho$ denote the composition
 \[\begin{array}{c}
\{{\rm Id}_{1}\}\times \mathbb{R}\times Spin(r)\\
\downarrow\\
\{{\rm Id}_{1}\}\times SO(2)\times Spin(r)\\
\downarrow\\
{\{{\rm Id}_{1}\}\times SO(2)\times Spin(r)\over \mathbb{Z}_2}
   \end{array}
\]
Thus,
\begin{eqnarray*}
\ker(\tilde\rho)
   &=& \left<({\rm Id}_{1},\pi, {\rm vol}_r)\right> ,
\\
\pi_1((\widehat{SO(1)}\times\widehat{SO(2)})\widehat{Spin(r)})
   &\cong& \mathbb{Z}. 
\end{eqnarray*}

\end{itemize}

\item[(4)] 
Recall that $m_1,m_2>0$,  $m_1\equiv  m_2\equiv 0\,\,(\mbox{mod $2$})$.
\begin{itemize}
\item If $m_1,m_2\geq 4$,
let $\tilde\rho$ denote the composition
 \[\begin{array}{c}
Spin(m_1)\times Spin(m_2)\times Spin(r)\\
\downarrow\\
SO(m_1)\times SO(m_2)\times Spin(r)\\
\downarrow\\
{SO(m_1)\times SO(m_2)\times Spin(r)\over \mathbb{Z}_2\oplus \mathbb{Z}_2}
   \end{array}
\]
Thus,
{\tiny
\begin{eqnarray*}
\ker(\tilde\rho)
&=&
   \left\{   \begin{array}{ll}
\left<(-1,1,1),(1,-1,1),(1,{\rm vol}_{m_2},{\rm vol}_r),({\rm vol}_{m_1},1,-{\rm vol}_r)\right> & \mbox{if $m_1\equiv m_2\equiv 0$ (mod 4),} \\
\left<(-1,1,1),({\rm vol}_{m_1},1,-{\rm vol}_r),(1,{\rm vol}_{m_2},{\rm vol}_r)\right> & \mbox{if $m_1\equiv 0$ (mod 4) and $m_2\equiv 2$ (mod 4),} \\
 \left<(1,-1,1),({\rm vol}_{m_1},1,-{\rm vol}_r),(1,{\rm vol}_{m_2},{\rm vol}_r)\right> & \mbox{if $m_2\equiv 2$ (mod 4) and $m_2\equiv 0$ (mod 4),} \\
\left<({\rm vol}_{m_1},1,-{\rm vol}_r),(1,{\rm vol}_{m_2},{\rm vol}_r)\right> & \mbox{if $m_1\equiv m_2\equiv 2$ (mod 4),}
                                              \end{array}
\right.\\ 
\pi_1((\widehat{SO(m_1)}\times \widehat{SO(m_2)})\widehat{Spin(r)})
   &\cong&
   \left\{   \begin{array}{ll}
\mathbb{Z}_2\oplus\mathbb{Z}_2\oplus\mathbb{Z}_2\oplus\mathbb{Z}_2  & \mbox{if $m_1\equiv m_2\equiv 0$ (mod 4),} \\
\mathbb{Z}_2\oplus\mathbb{Z}_2\oplus\mathbb{Z}_4 & \mbox{if $m_1\equiv 0$ (mod 4) and $m_2\equiv 2$ (mod 4),} \\
\mathbb{Z}_2\oplus\mathbb{Z}_4\oplus\mathbb{Z}_2 & \mbox{if $m_2\equiv 2$ (mod 4) and $m_2\equiv 0$ (mod 4),} \\
\mathbb{Z}_4\oplus\mathbb{Z}_4 & \mbox{if $m_1\equiv m_2\equiv 2$ (mod 4).}
                                              \end{array}
\right. 
\end{eqnarray*}
}

\item If $m_1\geq 4, m_2=2$,
let $\tilde\rho$ denote the composition
 \[\begin{array}{c}
Spin(m_1)\times \mathbb{R}\times Spin(r)\\
\downarrow\\
SO(m_1)\times SO(2)\times Spin(r)\\
\downarrow\\
{SO(m_1)\times SO(2)\times Spin(r)\over \mathbb{Z}_2\oplus \mathbb{Z}_2}
   \end{array}
\]
Thus,
\begin{eqnarray*}
\ker(\tilde\rho)
&=&
\left\{
\begin{array}{ll}
\left<(-1,0,1),({\rm vol}_{m_1},0,-{\rm vol}_r),(1,\pi,{\rm vol}_r)\right> & \mbox{if $m_1\equiv 0$ (mod 4),} \\
\left<({\rm vol}_{m_1},0,-{\rm vol}_r),(1,\pi,{\rm vol}_r)\right>&  \mbox{if $m_1\equiv 2$ (mod 4),}
\end{array}
\right. \\
\pi_1((\widehat{SO(m_1)}\times \widehat{SO(2)})\widehat{Spin(r)})
&\cong&
\left\{
\begin{array}{ll}
\mathbb{Z}_2\oplus \mathbb{Z}_2\oplus\mathbb{Z} & \mbox{if $m_1\equiv 0$ (mod 4),} \\
\mathbb{Z}_4 \oplus \mathbb{Z}&  \mbox{if $m_1\equiv 2$ (mod 4).}
\end{array}
\right.
\end{eqnarray*}

\item If $m_1=2,m_2\geq 4$,
let $\tilde\rho$ denote the composition
 \[\begin{array}{c}
\mathbb{R}\times Spin(m_2)\times Spin(r)\\
\downarrow\\
SO(2)\times SO(m_2)\times Spin(r)\\
\downarrow\\
{SO(2)\times SO(m_2)\times Spin(r)\over \mathbb{Z}_2\oplus \mathbb{Z}_2}
   \end{array}
\]
Thus,
\begin{eqnarray*}
\ker(\tilde\rho)
&=&
\left\{
\begin{array}{ll}
\left<(0,-1,1),(0,{\rm vol}_{m_2},{\rm vol}_r),(\pi,1,-{\rm vol}_r)\right> & \mbox{if $m_2\equiv 0$ (mod 4),} \\
\left<(0,{\rm vol}_{m_2},{\rm vol}_r),(\pi,1,-{\rm vol}_r)\right> &  \mbox{if $m_2\equiv 2$ (mod 4),}
\end{array}
\right. \\
\pi_1((\widehat{SO(2)}\times \widehat{SO(m_2)})\widehat{Spin(r)})
&\cong&
\left\{
\begin{array}{ll}
\mathbb{Z}_2\oplus \mathbb{Z}_2\oplus\mathbb{Z} & \mbox{if $m_2\equiv 0$ (mod 4),} \\
\mathbb{Z}_4 \oplus \mathbb{Z}&  \mbox{if $m_2\equiv 2$ (mod 4).}
\end{array}
\right. 
\end{eqnarray*}

\item If $m_1=m_2=2$,
let $\tilde\rho$ denote the composition
 \[\begin{array}{c}
\mathbb{R}\times \mathbb{R}\times Spin(r)\\
\downarrow\\
SO(2)\times SO(2)\times Spin(r)\\
\downarrow\\
{SO(2)\times SO(2)\times Spin(r)\over \mathbb{Z}_2\oplus \mathbb{Z}_2}
   \end{array}
\]
Thus,
\begin{eqnarray*}
\ker(\tilde\rho)
&=&
\left<(0,\pi,{\rm vol}_r), (\pi,0,-{\rm vol}_r)\right> ,
\\
\pi_1((\widehat{SO(2)}\times \widehat{SO(2)})\widehat{Spin(r)})
&\cong& \mathbb{Z}\oplus\mathbb{Z}.
\end{eqnarray*}

\end{itemize}

\item[(5)] 
Recall that $m_1>0, m_2=0$, $m_1\equiv 1 $ (mod 2).
\begin{itemize}
\item If $m_1\geq 3$, 
let $\tilde\rho$ denote the composition
\[\begin{array}{c}
Spin(m_1)\times Spin(r)\\
\downarrow\\
SO(m_1)\times Spin(r)\\
\downarrow\\
{SO(m_1)\times Spin(r)\over \mathbb{Z}_2}
  \end{array}
\]
Thus
\begin{eqnarray*}
\ker(\tilde\rho)&=& \left<(-1,1),(1,{\rm vol}_r)\right>,\\
\pi_1(\widehat{SO(m_1)}\widehat{Spin(r)^+})&\cong& \mathbb{Z}_2\oplus\mathbb{Z}_2. 
\end{eqnarray*}

\item If $m_1=1$, 
let $\tilde\rho$ denote
\[\begin{array}{c}
\{{\rm Id}_{1}\}\times Spin(r)\\
\downarrow\\
{\{{\rm Id}_{1}\}\times Spin(r)\over \mathbb{Z}_2}
  \end{array}
\]
Thus,
\begin{eqnarray*}
\ker(\tilde\rho)&=& \left<({\rm Id}_{1},{\rm vol}_r)\right>,\\
\pi_1(\widehat{SO(1)}\widehat{Spin(r)^+})&\cong& \mathbb{Z}_2. 
\end{eqnarray*}

\end{itemize}

\item[(6)] 
Recall that $m_1>0, m_2=0$, $m_1\equiv 0 $ (mod 2).
\begin{itemize}
\item If $m_1\geq 4$, 
let $\tilde\rho$ denote the composition
\[\begin{array}{c}
Spin(m_1)\times Spin(r)\\
\downarrow\\
SO(m_1)\times Spin(r)\\
\downarrow\\
{SO(m_1)\times Spin(r)\over \mathbb{Z}_2\oplus\mathbb{Z}_2}
  \end{array}
\]
Thus,
\begin{eqnarray*}
\ker(\tilde\rho)&=& \left\{\begin{array}{ll}
\left<(-1,1),({\rm vol}_{m_1},-1),(1,{\rm vol}_r)\right>& \mbox{if $m_1\equiv 0$ (mod 4),}\\
\left<({\rm vol}_{m_1},-1),(1,{\rm vol}_r)\right>& \mbox{if $m_1\equiv 2$ (mod 4),}
                                   \end{array}
 \right.\\
\pi_1(\widehat{SO(m_1)}\widehat{Spin(r)^+})&\cong& \left\{\begin{array}{ll}
 \mathbb{Z}_2\oplus\mathbb{Z}_2\oplus\mathbb{Z}_2& \mbox{if $m_1\equiv 0$ (mod 4),}\\
  \mathbb{Z}_2\oplus\mathbb{Z}_4 & \mbox{if $m_1\equiv 2$ (mod 4).}
                                   \end{array}
 \right. 
\end{eqnarray*}

\item If $m_1=2$, 
let $\tilde\rho$ denote the composition
\[\begin{array}{c}
\mathbb{R}\times Spin(r)\\
\downarrow\\
SO(2)\times Spin(r)\\
\downarrow\\
{SO(2)\times Spin(r)\over \mathbb{Z}_2\oplus\mathbb{Z}_2}
  \end{array}
\]
Thus,
\begin{eqnarray*}
\ker(\tilde\rho)&=& \left<(\pi,-1),(0,{\rm vol}_r)\right>,\\
\pi_1(\widehat{SO(2)}\widehat{Spin(r)^+})&\cong& \mathbb{Z}\oplus \mathbb{Z}_2. 
\end{eqnarray*}

\end{itemize}

\item[(7)] 
Recall that $m_1=0, m_2>0$, $m_2\equiv 1 $ (mod 2).
\begin{itemize}
\item If $m_2\geq 3$, 
let $\tilde\rho$ denote the composition
\[\begin{array}{c}
Spin(m_2)\times Spin(r)\\
\downarrow\\
SO(m_2)\times Spin(r)\\
\downarrow\\
{SO(m_2)\times Spin(r)\over \mathbb{Z}_2}
  \end{array}
\]
Thus,
\begin{eqnarray*}
\ker(\tilde\rho)&=& \left<(1,-{\rm vol}_r)\right>,\\
\pi_1(\widehat{SO(1)}\widehat{Spin(r)^-})&\cong& \mathbb{Z}_2. 
\end{eqnarray*}

\end{itemize}

\item[(8)] 
Recall that $m_1=0, m_2>0$, $m_2\equiv 0 $ (mod 2).
\begin{itemize}
\item If $m_2\geq 4$, 
let $\tilde\rho$ denote the composition
\[\begin{array}{c}
Spin(m_2)\times Spin(r)\\
\downarrow\\
SO(m_2)\times Spin(r)\\
\downarrow\\
{SO(m_2)\times Spin(r)\over \mathbb{Z}_2\oplus\mathbb{Z}_2}
  \end{array}
\]
Thus,
\begin{eqnarray*}
\ker(\tilde\rho)&=& \left\{\begin{array}{ll}
\left<(-1,1),({\rm vol}_{m_2},-1),(1,-{\rm vol}_r)\right>& \mbox{if $m_2\equiv 0$ (mod 4),}\\
\left<({\rm vol}_{m_2},-1),(1,-{\rm vol}_r)\right> & \mbox{if $m_2\equiv 2$ (mod 4),}
                                   \end{array}
 \right.\\
\pi_1(\widehat{SO(m_2)}\widehat{Spin(r)^-})&\cong& \left\{\begin{array}{ll}
\mathbb{Z}_2\oplus\mathbb{Z}_2\oplus\mathbb{Z}_2& \mbox{if $m_2\equiv 0$ (mod 4),}\\
 \mathbb{Z}_2\oplus\mathbb{Z}_4 & \mbox{if $m_2\equiv 2$ (mod 4).}
                                   \end{array}
 \right. 
\end{eqnarray*}

\item If $m_2=2$, 
let $\tilde\rho$ denote the composition
\[\begin{array}{c}
\mathbb{R}\times Spin(r)\\
\downarrow\\
SO(2)\times Spin(r)\\
\downarrow\\
{SO(2)\times Spin(r)\over \mathbb{Z}_2\oplus\mathbb{Z}_2}
  \end{array}
\]
Thus,
\begin{eqnarray*}
\ker(\tilde\rho)&=& \left<(\pi,-1),(0,-{\rm vol}_r)\right>,\\
\pi_1(\widehat{SO(2)}\widehat{Spin(r)^-})&\cong& \mathbb{Z}\oplus \mathbb{Z}_2. 
\end{eqnarray*}

\end{itemize}

\end{itemize}

\subsection{$r\equiv 2,6 \,\,({\rm mod}\,\, 8)$}

\underline{Complexification}.
The volume form ${\rm vol}_r=e_1\cdots e_r$ acts as a complex structure $J$ on $\tilde\Delta_r$ (cf. \cite{Arizmendi-Herrera}).
Therefore 
\begin{eqnarray*}
 \tilde\Delta_r \otimes \mathbb{C}
    &=& \{\psi - i J\psi \,| \, \psi\in \tilde\Delta_r\} \oplus \{\psi + i J\psi \,| \, \psi\in \tilde\Delta_r\}.
\end{eqnarray*}
Note that the involution $(-i)^{r/2}{\rm vol}_r\cdot$ acts on $\tilde\Delta_r\otimes\mathbb{C}$  as follows:
\begin{itemize}
 \item if $r\equiv 2\,\,({\rm mod}\,\, 8)$,
\begin{eqnarray*}
 (-i)^{r/2}{\rm vol}_r\cdot (\psi - i J\psi)
   &=& 
 (\psi - i J\psi) ,\\
 (-i)^{r/2}{\rm vol}_r\cdot (\psi + i J\psi)
   &=&
 -(\psi + i J\psi) ,
\end{eqnarray*}
i.e.
\begin{eqnarray*}
 \tilde\Delta_r \otimes \mathbb{C}
    &=& 
   \{\psi - i J\psi \,| \, \psi\in \tilde\Delta_r\} \oplus \{\psi + i J\psi \,| \, \psi\in \tilde\Delta_r\}\\
    &=&
   \Delta_r^+ \oplus \Delta_r^-.
\end{eqnarray*}

 \item if $r\equiv 6\,\,({\rm mod}\,\, 8)$,
\begin{eqnarray*}
 (-i)^{r/2}{\rm vol}_r\cdot (\psi - i J\psi)
   &=&
 -(\psi - i J\psi) ,\\
 (-i)^{r/2}{\rm vol}_r\cdot (\psi + i J\psi)
   &=&
 (\psi + i J\psi) ,
\end{eqnarray*}
i.e.
\begin{eqnarray*}
 \tilde\Delta_r \otimes \mathbb{C}
    &=& 
   \{\psi - i J\psi \,| \, \psi\in \tilde\Delta_r\} \oplus \{\psi + i J\psi \,| \, \psi\in \tilde\Delta_r\}\\
    &=&
   \Delta_r^- \oplus \Delta_r^+.
\end{eqnarray*}

\end{itemize}
In any case,
\[\tilde\Delta_r\otimes\mathbb{C} = \Delta_r.\]
The centralizer subalgebra of $\widehat{\mathfrak{spin}(r)}$ in $\mathfrak{so}(N)$ is
\[C_{\mathfrak{so}(N)}(\widehat{\mathfrak{spin}(r)})=\widehat{\mathfrak{u}(m)}
=\mathfrak{so}(m)\otimes {\rm Id}_{\tilde\Delta_r}\oplus \sym^2\mathbb{R}^m\otimes J ,\] 
where $N=d_rm$,  $\mathfrak{so}(m)$ and 
$\sym^2\mathbb{R}^m$ act on $\mathbb{R}^m$ as skew-symmetric and symmetric endomorphisms respectively.
Let
\[v\otimes (\psi - iJ\psi) \in(\mathbb{R}^m\otimes \tilde\Delta_r) \otimes \mathbb{C}\]
and 
\[A\otimes {\rm Id}_{\tilde\Delta_r}+ B\otimes J \in \mathfrak{so}(m)\otimes {\rm Id}_{\tilde\Delta_r}\oplus \sym^2\mathbb{R}^m\otimes J.\]
Now,
\begin{eqnarray*}
 (A\otimes {\rm Id}_{\tilde\Delta_r}+ B\otimes J)(v\otimes (\psi - iJ\psi))
   &=&
   Av\otimes \psi - iAv\otimes J\psi    + Bv\otimes J \psi - iBv\otimes JJ\psi \\
   &=&
   Av\otimes \psi  + iBv\otimes \psi +(-i)(i)Bv\otimes J \psi  - iAv\otimes J\psi  \\
   &=&
   (A+iB)v\otimes (\psi-iJ\psi), 
\end{eqnarray*}
where $A+iB\in\mathfrak{u}(m)$. 
Similarly, for  $v\otimes (\psi + iJ\psi)$,
\begin{eqnarray*}
 (A\otimes {\rm Id}_{\tilde\Delta_r}+ B\otimes J)(v\otimes (\psi + iJ\psi))
   &=&
   (A-iB)v\otimes (\psi +iJ\psi)      .
\end{eqnarray*}
Thus, 
\begin{equation}
(\mathbb{R}^m\otimes \tilde\Delta_r) \otimes \mathbb{C} = 
\left\{ \begin{array}{ll}
\mathbb{C}^m\otimes \Delta_r^+ \oplus \overline{\mathbb{C}^m}\otimes \Delta_r^- & \mbox{if $r\equiv 2\,\,({\rm mod}\,\, 8)$},\\
\mathbb{C}^m\otimes \Delta_r^- \oplus \overline{\mathbb{C}^m}\otimes \Delta_r^+ & \mbox{if $r\equiv 6\,\,({\rm mod}\,\, 8)$},
        \end{array}
\right. \label{eq: splitting complexification r=2,6}
\end{equation}
where $\mathbb{C}^m$ is the standard representation of $U(m)$. Therefore, we have a representation
\[U(m)\times Spin(r) \longrightarrow SO(N)\subset {\rm Aut}((\mathbb{R}^m\otimes \tilde\Delta_r) \otimes \mathbb{C}).\]

\noi\underline{Structure group}.
Since $\widehat{\mathfrak{u}(m)}$ and $\widehat{\mathfrak{spin}(r)}$ commute with each other,
we can take separately the exponentials of their elements within $\mathbb{C}(N)$.
With respect to \rf{eq: splitting complexification r=2,6}, an element
$A\otimes {\rm Id}_{\tilde\Delta_r}+ B\otimes J\in\mathfrak{u}(m)$ 
looks as follows
\[\left(
\begin{array}{ll}
(A+ iB)\otimes {\rm Id}_{\Delta_r^\pm}  & \\
 & (A- iB)\otimes {\rm Id}_{\Delta_r^\mp}
\end{array}
\right)\]
so that the exponentials form
\[\left\{\left(
\begin{array}{ll}
e^{A+ iB}\otimes {\rm Id}_{\Delta_r^\pm}  & \\
 & e^{A- iB}\otimes {\rm Id}_{\Delta_r^\mp}
\end{array}
\right)\,\,:\,\, A\in\mathfrak{so}(m), B\in \sym^2\mathbb{R}^m\right\} =: \widehat{U(m)}.
\]
With respect to \rf{eq: splitting complexification r=2,6}, an element ${\rm Id}_{m}\otimes \xi\in\widehat{\mathfrak{spin}(r)}$, 
looks as follows
\[\left(\begin{array}{ll}
{\rm Id}_{m}\otimes \kappa_{r^*}^+(\xi) & \\
 & {\rm Id}_{m}\otimes \kappa_{r^*}^-(\xi)
        \end{array}
\right)={\rm Id}_{m}\otimes \xi,\]
and its exponential  is
\[{\rm Id}_{m}\otimes e^\xi \in {\rm Id}_{m}\otimes \kappa(Spin(r)) =: \widehat{Spin(r)} \cong Spin(r),\]
since $Spin(r)$ is represented faithfully on $\Delta_r$.
The image of $U(m)\times Spin(r)$ in $SO(N)\subset {\rm Aut}((\mathbb{R}^m\otimes \tilde\Delta_r) \otimes \mathbb{C})$
under the aforementioned representation is 
\[\mathcal{N}_{SO(N)}^0(\widehat{Spin(r)}) =\widehat{U(m)}\widehat{Spin(r)},\] 
the subgroup of all possible products of elements of the two subgroups, i.e.
we have a map
\[U(m)\times Spin(r) \xrightarrow{\rho}  \widehat{U(m)}\widehat{Spin(r)}\subset SO(N).\]
Now we need to find $\ker(\rho)$ and identify $\widehat{U(m)}\widehat{Spin(r)}$ as a quotient
\[\widehat{U(m)}\widehat{Spin(r)}\cong {U(m)\times Spin(r)\over \ker(\rho)} .  \]
If there are elements 
$g\in U(m)$ and $h \in Spin(r)$ such that
\[\rho(g,h)={\rm Id}_{N},\]
then 
\[\widehat{Spin(r)}\ni\rho({\rm Id}_{m},h)=\rho(g,1)^{-1}\in  \widehat{U(m)}.\]
Since 
$\rho({\rm Id}_{m},h)$ commutes with every element of $\widehat{Spin(r)}$, it belongs to the center 
$Z(\widehat{Spin(r)})\cong Z(Spin(r))=\{ 1,-1, {\rm vol}_r,-{\rm vol}_r\} = \left<{\rm vol}_r\right>\cong \mathbb{Z}_4$. Recall that ${\rm vol}_r= e_1\cdots e_r$  
acts as $\mp i$ on $\Delta_r^\pm$ if $r\equiv 2 \,\,({\rm mod}\,\, 8)$,
and as $\pm i$ on $\Delta_r^\pm$ if $r\equiv 6 \,\,({\rm mod}\,\, 8)$, so that
it maps to
\[\mp \,(i{\rm Id}_{\Delta_r^+}\oplus (-i){\rm Id}_{\Delta_r^-})\]
in the complex $Spin(r)$ representation. 
Note that 
$({\rm Id}_{m}, {\rm vol}_r)$ maps to
\begin{eqnarray*}
 (-i){\rm Id}_{m}\otimes{\rm Id}_{\Delta_r^+}\oplus (i){\rm Id}_{m}\otimes {\rm Id}_{\Delta_r^-} 
&&\quad \mbox{if $r\equiv 2 \,\,({\rm mod}\,\, 8)$,}\\
(i){\rm Id}_{m}\otimes{\rm Id}_{\Delta_r^+}\oplus (-i){\rm Id}_{m}\otimes {\rm Id}_{\Delta_r^-}
&&\quad\mbox{if $r\equiv 6 \,\,({\rm mod}\,\, 8)$,}
\end{eqnarray*}
in $SO(N)$, 
and that  $((-i){\rm Id}_{\mathbb{C}^m},1)\in U(m)\times Spin(r)$  maps to such transformations in both cases. Thus, the elements of 
$U(m)\times Spin(r)$ mapping to ${\rm Id}_{N}$ are
\[\pm ({\rm Id}_{m},1),\quad\quad \pm(i{\rm Id}_{m}, -{\rm vol}_r),\]
which form a copy of $\mathbb{Z}_4$ and
\[\widehat{U(m)}\widehat{Spin(r)}\cong {U(m)\times Spin(r)\over \mathbb{Z}_4}.\]

\noi \underline{Fundamental group}.
Let
\begin{eqnarray*}
\mathbb{R}\times SU(m)\times Spin(r) &\xrightarrow{\tilde\rho}& \widehat{U(m)} \widehat{Spin(r)}\\ 
 (t, A, g) &\mapsto& \left(\begin{array}{ll}
 e^{it}A\otimes \kappa_r^\pm(g) & \\
 &  e^{-it}\overline{A}\otimes \kappa_r^\mp(g)
                           \end{array}
\right),
\end{eqnarray*}
Thus
\begin{eqnarray*}
\ker(\tilde\rho) &=& \left<\left( {2\pi\over m}, e^{-{2\pi i \over m}}{\rm Id}_{m}, 1\right),\left({\pi\over 2}, {\rm Id}_{m}, -{\rm vol}\right)\right>,
\\
\pi_1(\widehat{U(m)}\widehat(Spin(r)))
&\cong&\left\{
\begin{array}{ll}
\mathbb{Z}, & \mbox{if $(m,4)=1$,}\\
\mathbb{Z}\oplus\mathbb{Z}_2, & \mbox{if $(m,4)=2$,}\\
\mathbb{Z}\oplus\mathbb{Z}_4, & \mbox{if $(m,4)=4$.}
\end{array} 
\right. 
\end{eqnarray*}
Indeed, let
\[a:=\left( {2\pi\over m}, e^{-{2\pi i \over m}}{\rm Id}_{m}, 1\right),\quad b:=\left({\pi\over 2}, {\rm Id}_{m}, -{\rm vol}\right),\]
and note that (in multiplicative notation)
\[a^m = b^4.\]
Moreover,
\begin{itemize}
\item If $(m,4)=1$, there exist $t,m\in \mathbb{Z}$ coprime such that
\[tm+s4=1.\]
The element
\[b^ta^s\]
is such that
\begin{eqnarray*}
(b^ta^s)^m
 &=& b^{mt}(b^4)^s\\
 &=& b,\\
(b^ta^s)^4
  &=& (a^m)^ta^{4s}\\ 
  &=& a.
\end{eqnarray*}

\item If $(m,4)=2$, $m=4k+2$ and there exist two generators
\begin{eqnarray*}
 c&=& a^{-(2k+1)}b^2,\\
 d&=& ba^{-k},
\end{eqnarray*}
such that
\begin{eqnarray*}
 c^2&=&1,\\
 a&=& d^2c,\\
 b&=& d^{2k+1}c^k.
\end{eqnarray*}

\item If $(m,4)=4$, $m=4k$ and we have two generators 
\[a\quad\mbox{and}\quad c=a^{-k}b,\]
such that
\[ c^4 =1.\]
\end{itemize}

\subsection{$r\equiv 3,5 \,\,({\rm mod}\,\, 8)$}

\underline{Complexification}. In this case $\tilde\Delta_r$ admits three complex structures $I,J,K$, described explicitly in \cite{Arizmendi-Herrera}, 
which behave like quaternions and commute with  $\mathfrak{spin}(r)$. 
Let us consider the complexification of $\tilde\Delta_r$ and decompose as follows 
\[\tilde\Delta_r\otimes \mathbb{C} = \{\psi - iI\psi \,|\,   \psi \in \tilde\Delta_r \} 
\oplus \{ \psi + iI\psi\,|\,  \psi \in \tilde\Delta_r \},\]
where the first and second subspaces are the $+i$ and $-i$ eigenspaces of $I$ respectively.
Notice that 
\begin{eqnarray*}
  J(\psi \mp i I\psi) 
    &=&
   J\psi \mp i  JI\psi \\
    &=&
    J\psi \pm i IJ\psi ,
\end{eqnarray*}
i.e. $J$ interchanges the two subspaces and squares to $-{\rm Id}_{d_r}$.
For any $ \xi\in\mathfrak{spin}(r)$
\begin{eqnarray*}
  \xi(\psi \pm iI\psi)
    &=&
  \xi\psi \pm i\xi I\psi \\ 
    &=&
  \xi\psi \pm iI\xi\psi , 
\end{eqnarray*}
which means that the subspaces $\{ \psi - iI\psi \,|\,   \psi \in \tilde\Delta_r \}$ and  
$\{ \psi + iI\psi \,|\,  \psi \in \tilde\Delta_r \}$
are irreducible complex representations of $\mathfrak{spin}(r)$ of dimension $d_r/2$. 
Thus, they are isomorphic to $\Delta_r$ as $\mathfrak{spin}(r)$ representations and
\[\tilde\Delta_r \otimes \mathbb{C} \cong \Delta_r\oplus \Delta_r.\]
Now recall that the centralizer subalgebra of $\widehat{\mathfrak{spin}(r)}$ is
\begin{eqnarray*}
C_{\mathfrak{so}(N)}(\widehat{\mathfrak{spin}(r)})&=&\mathfrak{so}(m)\otimes {\rm Id}_{\tilde\Delta_r} \oplus \sym^2\mathbb{R}^m\otimes I
\oplus \sym^2\mathbb{R}^m\otimes J\oplus \sym^2\mathbb{R}^m\otimes K\\
&\cong& \mathfrak{so}(m)
\otimes {\rm Id}_{\tilde\Delta_r} \oplus \sym^2\mathbb{R}^m\otimes \mathfrak{sp}(1) \\
&\cong& \mathfrak{sp}(m),
\end{eqnarray*}
where $N=d_rm$.
Let us consider the complexification of $\mathbb{R}^m\otimes\tilde\Delta_r$ and decompose it 
\begin{equation*}
(\mathbb{R}^m\otimes\tilde\Delta_r)\otimes \mathbb{C} = \{ v\otimes (\psi - iI\psi) \,|\, v\in\mathbb{R}^m,  \psi \in \tilde\Delta_r \} 
\oplus \{ v\otimes (\psi + iI\psi) \,|\, v\in\mathbb{R}^m, \psi \in \tilde\Delta_r \}, 
\end{equation*}
where the first and second subspaces are the $+i$ and $-i$ eigenspaces of ${\rm Id}_{m}\otimes I$ respectively.
Notice that 
\begin{eqnarray*}
 ({\rm Id}_{m}\otimes J)(v\otimes (\psi \mp i{\rm Id}_{m}\otimes I\psi)) 
    &=&
   v\otimes J\psi \mp i v\otimes JI\psi \\
    &=&
   v\otimes J\psi \pm iv\otimes IJ\psi ,
\end{eqnarray*}
i.e. ${\rm Id}_{m}\otimes J$ interchanges the two subspaces and squares to $-{\rm Id}_{N}$.
For any ${\rm Id}_{m}\otimes \xi\in\widehat{\mathfrak{spin}(r)}$
\begin{eqnarray*}
 ({\rm Id}_{m}\otimes \xi)(v\otimes(\psi \pm iI\psi))
    &=&
   v\otimes (\xi\psi \pm i\xi I\psi) \\ 
    &=&
   v\otimes (\xi\psi \pm iI\xi\psi) , 
\end{eqnarray*}
which means that the subspaces $\{ v\otimes (\psi - iI\psi) \,|\, v\in\mathbb{R}^m,  \psi \in \tilde\Delta_r \}$ and  
$\{ v\otimes (\psi + iI\psi) \,|\, v\in\mathbb{R}^m, \psi \in \tilde\Delta_r \}$
are isomorphic to $\mathbb{C}^m\otimes\Delta_r$ as $\widehat{\mathfrak{spin}(r)}$ representations.

Now consider
\[A\otimes {\rm Id}_{\tilde\Delta_r} + B\otimes I + C\otimes J + D\otimes K
\in\mathfrak{so}(m)\otimes {\rm Id}_{\tilde\Delta_r} \oplus \sym^2\mathbb{R}^m\otimes I
\oplus \sym^2\mathbb{R}^m\otimes J\oplus \sym^2\mathbb{R}^m\otimes K = \mathfrak{sp}(m),\]
and 
\begin{eqnarray*}
 (A\otimes {\rm Id}_{\tilde\Delta_r} + B\otimes I + C\otimes J + D\otimes K) ( v\otimes (\psi + iI\psi))
    &=&
   Av\otimes (\psi + iI\psi) 
  +Bv\otimes (I\psi + i II\psi)\\
    &&
  +Cv\otimes (J\psi + i JI\psi)
  +Dv\otimes (K\psi + i KI\psi)\\
    &=&
   ((A-iB)\otimes {\rm Id}_{\tilde\Delta_r} + (C+iD)\otimes J)(v\otimes (\psi + iI\psi)) .
\end{eqnarray*}
Similarly,
\begin{eqnarray*}
 (A\otimes {\rm Id}_{\tilde\Delta_r} + B\otimes I + C\otimes J + D\otimes K) ( v\otimes (\psi - iI\psi))
    &=&
   ((A+iB)\otimes {\rm Id}_{\tilde\Delta_r} + (C-iD)\otimes J)(v\otimes (\psi - iI\psi)) .
\end{eqnarray*}
If $C=D=0$, the subalgebra
\[\widehat{\mathfrak{u}(m)_{I}}=\{A\otimes {\rm Id}_{\tilde\Delta_r} + B\otimes I 
\in\mathfrak{so}(m)\otimes {\rm Id}_{\tilde\Delta_r} \oplus \sym^2\mathbb{R}^m\otimes I\, |\, A\in \mathfrak{so}(m), B\in \sym^2\mathbb{R}^m\}\]
is represented as follows
\begin{eqnarray*}
\tilde\Delta_r\otimes\mathbb{C} 
   &=& 
  \mathbb{C}_I^m\otimes \Delta_r\oplus\overline{\mathbb{C}_I^m}\otimes \Delta_r,\\
   &=& 
  (\mathbb{C}_I^m\oplus\overline{\mathbb{C}_I^m})\otimes \Delta_r,
\end{eqnarray*}
where $\mathbb{C}_I^m$ and $\overline{\mathbb{C}_I^m}$ denote the standard representation of $\widehat{\mathfrak{u}(m)_I}$ and its conjugate respectively. 
Since ${\rm Id}_m\otimes J$ interchanges the two summands, squares to $-{\rm Id}_{N}$ and commutes with the
action of $\widehat{\mathfrak{spin}(r)}$,
we have the  standard complex representation of $\mathfrak{sp}(m)$ as a factor
\begin{equation}
\tilde\Delta_r\otimes\mathbb{C} = \mathbb{C}^{2m}\otimes \Delta_r.\label{eq: splitting complexification r=3,5} 
\end{equation}
Thus, we have a representation 
\[Sp(m)\times Spin(r) \longrightarrow SO(N) \subset {\rm Aut}(\mathbb{C}^{2m}\otimes \Delta_r).\]

\noi\underline{Structure group}.
Since $\widehat{\mathfrak{sp}(m)}$ and $\widehat{\mathfrak{spin}(r)}$ commute with each other,
we can take separately the exponentials of their elements within $\mathbb{C}(N)$.
By considering \rf{eq: splitting complexification r=3,5}, the exponential of an element 
$\Omega\otimes{\rm Id}_{\Delta_r}\in\mathfrak{sp}(m)\otimes {\rm Id}_{\Delta_r}=\widehat{\mathfrak{spin}(r)}$
is
\[ e^{\Omega}\otimes{\rm Id}_{\tilde\Delta_r} \in
\widehat{Sp(m)}= Sp(m)\otimes{\rm Id}_{\Delta_r}\cong Sp(m).
 \]
On the other hand, if ${\rm Id}_{2m}\otimes \xi\in\widehat{\mathfrak{spin}(r)}$, 
its exponential  is
\[{\rm Id}_{2m}\otimes e^\xi \in {\rm Id}_{2m}\otimes \kappa(Spin(r)) = \widehat{Spin(r)} \cong Spin(r),\]
since $Spin(r)$ is represented faithfully on $\Delta_r$.
The image of $Sp(m)\times Spin(r)$ in $SO(N)\subset {\rm Aut}(\mathbb{C}^{2m} \otimes \Delta_r)$
under the aforementioned representation is 
\[\mathcal{N}_{SO(N)}^0(\widehat{Spin(r)}) =\widehat{Sp(m)}\widehat{Spin(r)},\] 
the subgroup of all possible products of elements of the two subgroups, i.e.
we have a map
\[Sp(m)\times Spin(r) \xrightarrow{\rho}  \widehat{Sp(m)}\widehat{Spin(r)}\subset SO(N).\]
Now we need to find $\ker(\rho)$ and identify $\widehat{Sp(m)}\widehat{Spin(r)}$ as a quotient
\[\widehat{Sp(m)}\widehat{Spin(r)}\cong {Sp(m)\times Spin(r)\over \ker(\rho)} .  \]
If there are elements 
$g\in Sp(m)$ and $h \in Spin(r)$ such that
\[\rho(g,h)={\rm Id}_{N},\]
then 
\[\widehat{Spin(r)}\ni\rho({\rm Id}_{2m},h)=\rho(g,1)^{-1}\in \widehat{Sp(m)}.\]
Since
$\rho({\rm Id}_{2m},h)$ commutes with every element of $\widehat{Spin(r)}$, it belongs to its center 
$Z(\widehat{Spin(r)})\cong Z(Spin(r)) = \mathbb{Z}_2 =\{\pm 1\}$. Note that
$-1$ is mapped to $-{\rm Id}_{\Delta_r}$ under the $Spin(r)$ representation, and 
that $({\rm Id}_{2m},-1)$ maps to 
$-{\rm Id}_{2m}\otimes{\rm Id}_{\Delta_r}\in SO(N) $ under $\rho$.
Note that $-{\rm Id}_{2m}\otimes{\rm Id}_{\Delta_r}$ also belongs to $\widehat{Sp(m)}$
being the image of $(-{\rm Id}_{2m},1) \in Sp(m)\times Spin(r)$.
Thus,
\[\ker(\rho)=\{\pm ({\rm Id}_{2m},1) \}\cong \mathbb{Z}_2,\]
\[\widehat{Sp(m)}\widehat{Spin(r)} \cong {Sp(m)\times Spin(r)\over \mathbb{Z}_2} . \]

\noi \underline{Fundamental group}. Clearly, 
\[\pi_1(\widehat{Sp(m)}\widehat{Spin(r)})=\mathbb{Z}_2 . \]

\subsection{$r\equiv 4 \,\,({\rm mod}\,\, 8)$}

Recall from \cite{Arizmendi-Herrera} that
\[\tilde\Delta_r^\pm = {1\over 2}(1\pm e_1\cdots e_r)\tilde\Delta_{r+3}.\]
In this case, $\tilde\Delta_r^\pm$ admits three complex structure $I^\pm,J^\pm$ and $K^\pm$ 
induced by Clifford multiplication with the elements 
$\frac12(1\pm e_{1}\dots e_{r})e_{r+1}e_{r+2}$, $\frac12(1\pm e_{1}\dots e_{r})e_{r+1}e_{r+3}$
and $\frac12(1\pm e_{1}\dots e_{r})e_{r+2}e_{r+3}$, respectively.
Just as in the previous case, 
\[\tilde\Delta_r^+\otimes \mathbb{C} = \{ \psi - iI^+\psi \,|\, \psi \in \tilde\Delta_r^+ \} \oplus \{ \psi + iI^+\psi \,|\, \psi \in \tilde\Delta_r^+ \}, \]
and both summands are isomorphic to $\Delta_r^-$.
Indeed, if
\begin{eqnarray*}
 \psi &=& {1\over 2}(1+ e_1\cdots e_r)\cdot \phi\in\tilde\Delta_r^+, 
\end{eqnarray*}
then,
\begin{eqnarray*}
(-i)^{r/2}(e_1\cdots e_r)\cdot (\psi \pm iI^+\psi )
   &=& - (\psi \pm iI^+\psi ),
\end{eqnarray*}
i.e. $\psi \pm iI^+\psi \in\Delta_r^-$. In other words,
\[\tilde\Delta_r^+\otimes\mathbb{C} = \Delta_r^-.\]
Similarly,
\[\tilde\Delta_r^-\otimes\mathbb{C} = \Delta_r^+.\]
The rest of the proof proceeds as in the previous case, 
\[(\mathbb{R}^{m_1}\otimes \tilde\Delta_r^+\oplus \mathbb{R}^{m_2}\tilde\Delta_r^-)\otimes \mathbb{C}=
\mathbb{C}^{2m_1}\otimes \Delta_r^-\oplus \mathbb{C}^{2m_2}\Delta_r^+,\]
and we have a representation
\[Sp(m_1)\times Sp(m_2)\times Spin(r) \longrightarrow SO(N)\subset {\rm Aut}(\mathbb{C}^N),\]
where $N=d_r(m_1+m_2)$.

\noi\underline{Structure group}.
Since $\widehat{\mathfrak{sp}(m_1)}\oplus\widehat{\mathfrak{sp}(m_2)}$ and $\widehat{\mathfrak{spin}(r)}$ commute with each other,
we can take separately the exponentials of their elements within $\mathbb{C}(N)$.
The exponential of 
$\Omega_1\otimes {\rm Id}_{\Delta_r^-}\oplus \Omega_2\otimes {\rm Id}_{\Delta_r^+}\in \widehat{\mathfrak{sp}(m_1)}\oplus\widehat{\mathfrak{sp}(m_2)}$ gives 
\begin{eqnarray*}
\left(
\begin{array}{cc}
e^{\Omega_1}\otimes{\rm Id}_{\Delta_r^-} & \\
 & e^{\Omega_2}\otimes{\rm Id}_{\Delta_r^+}
\end{array}
\right) 
 &\in&
\left(
\begin{array}{cc}
\widehat{Sp(m_1)} & \\
 & \widehat{Sp(m_2)}
\end{array}
\right)\\
&=&
\left(
\begin{array}{cc}
Sp(m_1)\otimes{\rm Id}_{\Delta_r^-} & \\
 & Sp(m_2)\otimes{\rm Id}_{\Delta_r^+}
\end{array}
\right)\\
 &\cong& Sp(m_1)\times Sp(m_2).
\end{eqnarray*}
On the other hand, if ${\rm Id}_{2m_1}\otimes \xi^-\oplus {\rm Id}_{2m_2}\otimes \xi^+\in\widehat{\mathfrak{spin}(r)}$, 
its exponential  is
\begin{eqnarray*}
\left(
\begin{array}{cc}
{\rm Id}_{2m_1}\otimes e^{\xi^-} & \\
 & {\rm Id}_{2m_2}\otimes e^{\xi^+}
\end{array}
\right)
  &\in& 
  \left(
\begin{array}{cc}
{\rm Id}_{2m_1}\otimes \kappa^-(Spin(r)) & \\
 & {\rm Id}_{2m_2}\otimes \kappa^+(Spin(r))
\end{array}\right)\\
&=& 
\left\{
\begin{array}{ll}
\widehat{Spin(r)}\cong Spin(r) & \mbox{if $m_1>0$ and $m_2>0$,}\\
\widehat{Spin(r)^-}\cong \kappa^-(Spin(r))& \mbox{if $m_1>0$ and $m_2=0$,}\\
\widehat{Spin(r)^+}\cong \kappa^+(Spin(r)) & \mbox{if $m_1=0$ and $m_2>0$,}
\end{array}
\right.
\end{eqnarray*}
where the first case is faithful and the last two are not, with 
\begin{eqnarray*}
 \widehat{Spin(r)^\pm}&\cong& \kappa^\pm(Spin(r)) \cong {Spin(r)\over \{1, \mp{\rm vol}_r\} }, \quad \mbox {if $r>4$}\\
 \widehat{Spin(r)^\pm}&\cong& \kappa^\pm(Spin(r))\cong  Spin(3) , \quad\quad \mbox {if $r=4$.}
\end{eqnarray*}
The images of $Sp(m_1)\times Sp(m_2)\times Spin(r)$ in $SO(N)\subset {\rm Aut}(\mathbb{C}^{N} )$
under the aforementioned representations are 
\begin{eqnarray*}
\mathcal{N}_{SO(N)}^0(\widehat{Spin(r)}) &=&(\widehat{Sp(m_1)}\times \widehat{Sp(m_2)})\widehat{Spin(r)},\\
\mathcal{N}_{SO(N)}^0(\widehat{Spin(r)^-}) &=&\widehat{Sp(m_1)}\widehat{Spin(r)^-}, \\ 
\mathcal{N}_{SO(N)}^0(\widehat{Spin(r)^+}) &=&\widehat{Sp(m_2)}\widehat{Spin(r)^+}, 
\end{eqnarray*}
respectively,
i.e. we have maps
\begin{eqnarray*}
Sp(m_1)\times Sp(m_2)\times Spin(r) &\xrightarrow{\rho}&  (\widehat{Sp(m_1)}\times \widehat{Sp(m_2)})\widehat{Spin(r)}\subset SO(N),\\ 
Sp(m_1)\times Spin(r) &\xrightarrow{\rho}&  \widehat{Sp(m_1)}\widehat{Spin(r)^-}\subset SO(N),\\ 
Sp(m_2)\times Spin(r) &\xrightarrow{\rho} & \widehat{Sp(m_2)}\widehat{Spin(r)^+}\subset SO(N). 
\end{eqnarray*}
Now we need to find $\ker(\rho)$ en each case to identify the relevant group as a quotient.
\begin{itemize}
 \item Case $m_1, m_2>0$. If there are elements 
$g_i\in Sp(m_i)$ and $h \in Spin(r)$ such that
\[\rho(g_1,g_2,h)={\rm Id}_{N},\]
then 
\[\widehat{Spin(r)}\ni\rho({\rm Id}_{2m_1},{\rm Id}_{2m_2},h)=\rho(g_1,g_2,1)^{-1}\in \widehat{Sp(m_1)}\times \widehat{Sp(m_2)}.\]
Since $\rho({\rm Id}_{2m_1},{\rm Id}_{2m_2},h)$ commutes with every element of $\widehat{Spin(r)}$, it belongs to its center 
$Z(\widehat{Spin(r)})\cong Z(Spin(r))  =\{1,-1, {\rm vol}_r, -{\rm vol}_r\} \cong \mathbb{Z}_2 \oplus \mathbb{Z}_2 $.
The element $-1$ is mapped to $-{\rm Id}_{\Delta_r^\pm}$ in the $Spin(r)$ representations $\Delta_r^\pm$, and $({\rm Id}_{2m_1},{\rm Id}_{2m_2},-1)$ 
is mapped to 
$-({\rm Id}_{2m_1}\otimes{\rm Id}_{\Delta_r^-}\oplus {\rm Id}_{2m_2}\otimes{\rm Id}_{\Delta_r^+})\in SO(N) $.
The element ${\rm vol}_r$ is mapped to $\mp {\rm Id}_{\Delta_r^\pm}$ in the $Spin(r)$ representations $\Delta_r^\pm$, 
and $({\rm Id}_{2m_1},{\rm Id}_{2m_2},{\rm vol}_r)$ is mapped to
$({\rm Id}_{2m_1}\otimes{\rm Id}_{\Delta_r^-}\oplus (-1){\rm Id}_{2m_2}\otimes{\rm Id}_{\Delta_r^+})\in SO(d_r(m_1+m_2)) $.
In this case,
$-({\rm Id}_{2m_1}\otimes{\rm Id}_{\Delta_r^-}\oplus {\rm Id}_{2m_2}\otimes{\rm Id}_{\Delta_r^+})$  and
$({\rm Id}_{2m_1}\otimes{\rm Id}_{\Delta_r^-}\oplus (-1){\rm Id}_{2m_2}\otimes{\rm Id}_{\Delta_r^+})$
belong to $\widehat{Sp(m_1)}\times\widehat{Sp(m_2)}$.
Thus,
 \begin{eqnarray*}
 \ker(\rho)&=&\{({\rm Id}_{2m_1},{\rm Id}_{2m_2},1),
     (-{\rm Id}_{2m_1}, -{\rm Id}_{2m_2},-1),\\
     &&({\rm Id}_{2m_1},-{\rm Id}_{2m_2},{\rm vol}_r),
     (-{\rm Id}_{2m_1},{\rm Id}_{2m_2},-{\rm vol}_r)\},\\
(\widehat{Sp(m_1)}\times \widehat{Sp(m_2)})\widehat{Spin(r)}&\cong& {Sp(m_1)\times Sp(m_2)\times Spin(r)\over \mathbb{Z}_2\oplus \mathbb{Z}_2} . 
 \end{eqnarray*}

\item Case $m_1>0$, $m_2=0$. If there are elements 
$g_1\in Sp(m_1)$ and $h \in Spin(r)$ such that
\[\rho(g_1,h)={\rm Id}_{N},\]
then 
\[\rho({\rm Id}_{2m_1},h)=\rho(g_1,1)^{-1}\in \widehat{Sp(m_1)}\]
and
\[\rho({\rm Id}_{2m_1},h)\in \widehat{Spin(r)^-}\cap \widehat{Sp(m_1)}.\]
Since $\rho({\rm Id}_{2m_1},h)$ commutes with every element of $\widehat{Spin(r)^-}$, it belongs to its center 
{\footnotesize
\[Z(\widehat{Spin(r)^-})\cong 
\left\{\begin{array}{ll}
Z(\kappa^-(Spin(r)))= Z(Spin(r)/\{1,{\rm vol}_r\})  = \{1,-1,{\rm vol}_r,-{\rm vol}_r\}/\{1, {\rm vol}_r\} \cong \{1,-1\} \cong \mathbb{Z}_2 & \mbox{if $r>4$.}\\
Z(\kappa^-(Spin(r)))= Z(\{1\}\times Spin(3))  = \{(1,1),(1,-1) \}  \cong \mathbb{Z}_2 &  \mbox{if $r=4$.} 
       \end{array}  
\right.
\]
}
\begin{itemize}
\item
If $r>4$, the element $-1$ is mapped to $-{\rm Id}_{\Delta_r^-}$ in the $Spin(r)$ representation $\Delta_r^-$, and $({\rm Id}_{2m_1},-1)$ 
is mapped to 
$-({\rm Id}_{2m_1}\otimes{\rm Id}_{\Delta_r^-})\in SO(N) $.
In this case,
$-{\rm Id}_{2m_1}\otimes{\rm Id}_{\Delta_r^-}$  
belongs to $\widehat{Sp(m_1)}$.
Thus,
 \begin{eqnarray*}
 \ker(\rho)&=&\{({\rm Id}_{2m_1},1),({\rm Id}_{2m_1},{\rm vol}_r),
     (-{\rm Id}_{2m_1},-1),(-{\rm Id}_{2m_1},-{\rm vol}_r)\},\\
\widehat{Sp(m_1)}\widehat{Spin(r)^-}&\cong& {Sp(m_1)\times Spin(r)\over  \mathbb{Z}_2\oplus \mathbb{Z}_2} . 
 \end{eqnarray*}

\item
If $r=4$, the element $(1,-1)\in \{1\}\times Spin(3)$ is mapped to $-{\rm Id}_{\Delta_r^-}$ in the $Spin(r)$ representation $\Delta_r^-$, and $({\rm Id}_{2m_1},(1,-1))$ 
is mapped to 
$-({\rm Id}_{2m_1}\otimes{\rm Id}_{\Delta_r^-})\in SO(N) $.
In this case,
$-{\rm Id}_{2m_1}\otimes{\rm Id}_{\Delta_r^-}$  
belongs to $\widehat{Sp(m_1)}$.
Thus,
 \begin{eqnarray*}
 \ker(\rho)&=& \{({\rm Id}_{2m_1},(1,1)),(-{\rm Id}_{2m_1},(1,-1))\}\times (Spin(3)\times \{1\}),\\
\widehat{Sp(m_1)}\widehat{Spin(4)^-}&\cong& {Sp(m_1)\times (Spin(3)\times Spin(3))\over   
\{({\rm Id}_{2m_1},(1,1)),(-{\rm Id}_{2m_1},(1,-1))\}\times (Spin(3)\times \{1\})}  \cong {Sp(m_1)\times Spin(3)\over \mathbb{Z}_2} . 
 \end{eqnarray*}
Note that $(1,-1)\in Spin(3)\times Spin(3)$ corresponds $-{\rm vol}_4\in Spin(4)$.
 
 \end{itemize}
 
\item Case $m_1=0$, $m_2>0$. 
If there are elements 
$g_2\in Sp(m_2)$ and $h \in Spin(r)$ such that
\[\rho(g_2,h)={\rm Id}_{N},\]
then 
\[\rho({\rm Id}_{2m_2},h)=\rho(g_2,1)^{-1}\in \widehat{Sp(m_2)}\]
and
\[\rho({\rm Id}_{2m_2},h)\in \widehat{Spin(r)^+}\cap \widehat{Sp(m_2)}.\]
Since $\rho({\rm Id}_{2m_2},h)$ commutes with every element of $\widehat{Spin(r)^+}$, it belongs to its center 
{\footnotesize
\[Z(\widehat{Spin(r)^-})\cong 
\left\{\begin{array}{ll}
Z(\kappa^+(Spin(r)))= Z(Spin(r)/\{1,-{\rm vol}_r\})  = \{1,-1,{\rm vol}_r,-{\rm vol}_r\}/\{1, -{\rm vol}_r\} \cong \{1,-1\} \cong \mathbb{Z}_2 & \mbox{if $r>4$.}\\
Z(\kappa^+(Spin(r)))= Z(Spin(3)\times \{1\})  = \{(1,1),(-1,1) \}  \cong \mathbb{Z}_2 &  \mbox{if $r=4$.} 
       \end{array}  
\right.
\]
}
\begin{itemize}
\item
The element $-1$ is mapped to $-{\rm Id}_{\Delta_r^+}$ in the $Spin(r)$ representations $\Delta_r^+$, and $({\rm Id}_{2m_1},-1)$ 
is mapped to 
$-({\rm Id}_{2m_2}\otimes{\rm Id}_{\Delta_r^+})\in SO(N) $.
In this case,
$-{\rm Id}_{2m_2}\otimes{\rm Id}_{\Delta_r^-}$  
belongs to $\widehat{Sp(m_2)}$.
Thus,
 \begin{eqnarray*}
 \ker(\rho)&=&\{({\rm Id}_{2m_2},1),({\rm Id}_{2m_2},-{\rm vol}_r),
     (-{\rm Id}_{2m_2},-1),(-{\rm Id}_{2m_2},{\rm vol}_r)
     \},\\
\widehat{Sp(m_2)}\widehat{Spin(r)^+}&\cong& {Sp(m_2)\times Spin(r)\over  \mathbb{Z}_2\oplus\mathbb{Z}_2} . 
 \end{eqnarray*}
 
\item
If $r=4$, the element $(-1,1)\in Spin(3)\times \{1\}$ is mapped to $-{\rm Id}_{\Delta_r^+}$ in the $Spin(r)$ representation $\Delta_r^+$, and $({\rm Id}_{2m_2},(-1,1))$ 
is mapped to 
$-({\rm Id}_{2m_2}\otimes{\rm Id}_{\Delta_r^+})\in SO(N) $.
In this case,
$-{\rm Id}_{2m_2}\otimes{\rm Id}_{\Delta_r^-}$  
belongs to $\widehat{Sp(m_2)}$.
Thus,
 \begin{eqnarray*}
 \ker(\rho)&=& \{({\rm Id}_{2m_2},(1,1)),(-{\rm Id}_{2m_2},(-1,1))\}\times (Spin(3)\times \{1\}),\\
\widehat{Sp(m_2)}\widehat{Spin(4)^+}&\cong& {Sp(m_2)\times (Spin(3)\times Spin(3))\over   
\{({\rm Id}_{2m_2},(1,1)),(-{\rm Id}_{2m_2},(-1,1))\}\times (\{1\}\times Spin(3))}  \cong {Sp(m_2)\times Spin(3)\over \mathbb{Z}_2} . 
 \end{eqnarray*}
Note that $(-1,1)\in Spin(3)\times Spin(3)$ corresponds ${\rm vol}_4\in Spin(4)$.

\end{itemize}
 
\end{itemize}

\noi \underline{Fundamental group}. Clearly,
\begin{eqnarray*}
\pi_1((\widehat{Sp(m_1)}\times\widehat{Sp(m_2)})\widehat{Spin(r)})&=&\mathbb{Z}_2\oplus \mathbb{Z}_2,\\
\pi_1(\widehat{Sp(m_1)}\widehat{Spin(r)^-})&\cong&
\left\{
\begin{array}{ll}
\mathbb{Z}_2\oplus \mathbb{Z}_2, & \mbox{if $m_1>0$, $m_2=0$, $r>4$,}\\
\mathbb{Z}_2, & \mbox{if $m_1>0$, $m_2=0$, $r=4$,}
\end{array}
\right.\\ 
\pi_1(\widehat{Sp(m_2)}\widehat{Spin(r)^+})&\cong& \left\{
\begin{array}{ll}
\mathbb{Z}_2\oplus \mathbb{Z}_2, & \mbox{if $m_1=0$, $m_2>0$, $r>4$,}\\
\mathbb{Z}_2 & \mbox{if $m_1=0$, $m_2>0$, $r=4$.}
\end{array}
\right.
\end{eqnarray*}

\qd

\vspace{.2in}

Thus, we have proved 
the following three theorems.

\begin{theo}\label{theo: complexifications}
The complexification of a real representation $\mathbb{R}^N$ of $Cl_r^0$ without trivial summands decomposes as follows
\[  
\begin{array}{|c|c|}
 \hline
  r \mbox{\ {\rm (mod 8)} }  &\mathbb{R}^N\otimes \mathbb{C}
\rule{0pt}{3ex}\\
\hline
 0 & \mathbb{C}^{m_1}\otimes \Delta_r^+ \oplus \mathbb{C}^{m_2}\otimes \Delta_r^-\tstrut\\
 \hline
 1,7 & \mathbb{C}^m\otimes \Delta_r\tstrut\\
\hline
 2 & \mathbb{C}^m\otimes \Delta_r^+\oplus\overline{\mathbb{C}^m}\otimes\Delta_r^-\tstrut\\
\hline
 6 & \overline{\mathbb{C}^m}\otimes \Delta_r^+\oplus \mathbb{C}^m\otimes\Delta_r^-\tstrut\\
\hline
 3,5 & \mathbb{C}^{2m}\otimes \Delta_r\tstrut\\
\hline
 4 & \mathbb{C}^{2m_2}\otimes \Delta_r^+ \oplus \mathbb{C}^{2m_1}\otimes \Delta_r^-\tstrut\\
\hline
  \end{array}
\]
where the different $\mathbb{C}^s$ denote the corresponding standard complex representations of the classical Lie algebras 
$\mathfrak{so}(s), \mathfrak{u}(s)$ or $\mathfrak{sp}(s/2)$.
\end{theo}
\qd

\begin{theo}\label{theo: quotient groups} The connected components of the identity $\mathcal{N}^0_{SO(N)}(\widehat{Spin(r)})$ 
of the normalizers $\mathcal{N}_{SO(N)}(\widehat{Spin(r)})$ are isomorphic to the following groups:
\begin{itemize} 
 \item If $r\equiv 1,7 \,\,({\rm mod}\,\, 8)$, $N=d_rm$ and
\[\mathcal{N}^0_{SO(N)}(\widehat{Spin(r)})
\cong \left\{
 \begin{array}{ll}
{SO(m)\times Spin(r)\over \mathbb{Z}_2}  ,&   \mbox{if $m$ is even,}\\ 
SO(m)\times Spin(r), & \mbox{if $m$ is odd,} 
\end{array}
 \right.\]

 \item If $r\equiv 0 \,\,({\rm mod}\,\, 8)$, $N=d_r(m_1+m_2)$ and 
 \begin{eqnarray*}
\mathcal{N}^0_{SO(N)}(\widehat{Spin(r)})
&\cong&
\left\{
\begin{array}{ll}
 {SO(m_1)\times SO(m_2)\times Spin(r)}, & \mbox{if $m_1>0, m_2>0$, $m_1\equiv m_2\equiv 1\,\,(\mbox{\rm mod $2$})$,}\\
 {SO(m_1)\times SO(m_2)\times Spin(r)\over \mathbb{Z}_2}, & \mbox{if $m_1>0, m_2>0$, $m_1+m_2\equiv 1\,\,(\mbox{\rm mod $2$})$,}\\
 {SO(m_1)\times SO(m_2)\times Spin(r)\over \mathbb{Z}_2\oplus \mathbb{Z}_2}, & \mbox{if $m_1>0, m_2>0$, $m_1\equiv  m_2\equiv 0\,\,(\mbox{\rm mod $2$})$,}
 \end{array}
\right.
\\
\mathcal{N}^0_{SO(N)}(\widehat{Spin(r)^+})
&\cong&
\left\{
\begin{array}{ll}
 {SO(m_1)\times Spin(r)\over \mathbb{Z}_2}, & \mbox{if $m_1>0, m_2=0$, $m_1\equiv 1\,\,(\mbox{\rm mod $2$})$,}\\
 {SO(m_1)\times Spin(r)\over \mathbb{Z}_2\oplus\mathbb{Z}_2}, & \mbox{if $m_1>0, m_2=0$, $m_1\equiv 0\,\,(\mbox{\rm mod $2$})$,}\\
 \end{array}
\right.\\
\mathcal{N}^0_{SO(N)}(\widehat{Spin(r)^-})
&\cong&
\left\{
\begin{array}{ll}
 {SO(m_2)\times Spin(r)\over \mathbb{Z}_2}, & \mbox{if $m_1=0, m_2>0$, $m_2\equiv 1\,\,(\mbox{\rm mod $2$})$,}\\
 {SO(m_2)\times Spin(r)\over \mathbb{Z}_2\oplus\mathbb{Z}_2}, & \mbox{if $m_1=0, m_2>0$, $m_2\equiv 0\,\,(\mbox{\rm mod $2$})$.}\\
 \end{array}
\right.
\end{eqnarray*}

\item If $r\equiv 2,6 \,\,({\rm mod}\,\, 8)$, $N=d_rm$ and
\[\mathcal{N}^0_{SO(N)}(\widehat{Spin(r)})
\cong {U(m)\times Spin(r)\over \mathbb{Z}_4}.\]
 
 \item If $r\equiv 3,5 \,\,({\rm mod}\,\, 8)$, $N=d_rm$ and
 \[\mathcal{N}^0_{SO(N)}(\widehat{Spin(r)})
 \cong {Sp(m)\times Spin(r)\over \mathbb{Z}_2} . \]

 \item If $r\equiv 4 \,\,({\rm mod}\,\, 8)$, $N=d_r(m_1+m_2)$ and
 \begin{eqnarray*}
 \mathcal{N}^0_{SO(N)}(\widehat{Spin(r)}) 
 &\cong&{Sp(m_1)\times Sp(m_2)\times Spin(r)\over \mathbb{Z}_2\oplus \mathbb{Z}_2} ,\quad \mbox{if $m_1>0$, $m_2>0$,}\\  
 \mathcal{N}^0_{SO(N)}(\widehat{Spin(r)^-}) 
 &\cong&\left\{
 \begin{array}{ll}
{Sp(m_1)\times Spin(r)\over \mathbb{Z}_2\oplus \mathbb{Z}_2}  & \mbox{if $m_1>0$, $m_2=0$, $r>4$,}\\
{Sp(m_1)\times Spin(3)\over \mathbb{Z}_2} & \mbox{if $m_1>0$, $m_2=0$, $r=4$,}
 \end{array}
 \right.\\  
 \mathcal{N}^0_{SO(N)}(\widehat{Spin(r)^+}) 
 &\cong& \left\{
\begin{array}{ll}
{Sp(m_2)\times Spin(r)\over \mathbb{Z}_2\oplus \mathbb{Z}_2} & \mbox{if $m_1=0$, $m_2>0$, $r>4$,}\\
{Sp(m_2)\times Spin(3)\over \mathbb{Z}_2} & \mbox{if $m_1=0$, $m_2>0$, $r=4$.}
\end{array}
\right. \\
 \end{eqnarray*}

\end{itemize}
 
\end{theo}
\qd

\begin{theo} \label{theo: fundamental groups}
The fundamental group of the connected components of the identity of the normalizers $\mathcal{N}^0_{SO(N)}(\widehat{Spin(r)})$ are the following.
\begin{itemize} 
 \item If $r\equiv 1,7 \,\,({\rm mod}\,\, 8)$, $N=d_rm$ and
\[
\pi_1(\mathcal{N}_{SO(N)}^0(\widehat{Spin(r)}))
\cong\left\{
 \begin{array}{ll}
\mathbb{Z}_2\oplus \mathbb{Z}_2, & \mbox{if $m\geq 4$, $m\equiv 0$ {\rm (mod 4)},} \\
\mathbb{Z}_4, & \mbox{if $m\geq 4$, $m\equiv 2$ {\rm (mod 4)},}\\
\mathbb{Z}_2, & \mbox{if $m> 1$ and odd,}\\ 
\{1\}, & \mbox{if $m= 1$,}\\ 
\mathbb{Z}, & \mbox{if $m=2$.} 
 \end{array}
\right.\]
 \item If $r\equiv 0 \,\,({\rm mod}\,\, 8)$, $N=d_r(m_1+m_2)$ and either $\pi_1(\mathcal{N}_{SO(N)}^0(\widehat{Spin(r)}))$ or $\pi_1(\mathcal{N}_{SO(N)}^0(\widehat{Spin(r)^+}))$ or 
 $\pi_1(\mathcal{N}_{SO(N)}^0(\widehat{Spin(r)^-}))$
are isomorphic to 
\vspace{-.15in}
\begin{center}
 {\rm 
 \begin{tabular}{|c|c|c|c|c|c|c|}\hline
\backslashbox{$m_1$}{$m_2$} & 0 & 1 & 2 & 1 (mod 2) & 2 (mod 4) & 0 (mod 4)\\\hline
0 &  & $\mathbb{Z}_2$ & $\mathbb{Z}\oplus\mathbb{Z}_2$ & $\mathbb{Z}_2\oplus\mathbb{Z}_2$ & $\mathbb{Z}_2\oplus\mathbb{Z}_4$ & $\mathbb{Z}_2\oplus\mathbb{Z}_2\oplus\mathbb{Z}_2$\\\hline
1 & $\mathbb{Z}_2$ & $\{1\}$ & $\mathbb{Z}$ & $\mathbb{Z}_2$ & $\mathbb{Z}_4$ & $\mathbb{Z}_2\oplus\mathbb{Z}_2$\\\hline
2 & $\mathbb{Z}\oplus\mathbb{Z}_2$ & $\mathbb{Z}$ & $\mathbb{Z}\oplus\mathbb{Z}$ & $\mathbb{Z}\oplus\mathbb{Z}_2$ & $\mathbb{Z}\oplus\mathbb{Z}_4$ & $\mathbb{Z}\oplus\mathbb{Z}_2\oplus\mathbb{Z}_2$\\\hline
1 (mod 2) & $\mathbb{Z}_2\oplus\mathbb{Z}_2$ & $\mathbb{Z}_2$ & $\mathbb{Z}\oplus\mathbb{Z}_2$ & $\mathbb{Z}_2\oplus\mathbb{Z}_2$ & $\mathbb{Z}_2\oplus\mathbb{Z}_4$ & $\mathbb{Z}_2\oplus\mathbb{Z}_2\oplus\mathbb{Z}_2$\\\hline
2 (mod 4) & $\mathbb{Z}_2\oplus\mathbb{Z}_4$ & $\mathbb{Z}_4$ & $\mathbb{Z}\oplus\mathbb{Z}_4$ & $\mathbb{Z}_2\oplus\mathbb{Z}_4$ & $\mathbb{Z}_4\oplus\mathbb{Z}_4$ & $\mathbb{Z}_2\oplus\mathbb{Z}_2\oplus\mathbb{Z}_4$\\\hline
0 (mod 4) & $\mathbb{Z}_2\oplus\mathbb{Z}_2\oplus\mathbb{Z}_2$ & $\mathbb{Z}_2\oplus\mathbb{Z}_2$ & $\mathbb{Z}\oplus\mathbb{Z}_2\oplus\mathbb{Z}_2$ & $\mathbb{Z}_2\oplus\mathbb{Z}_2\oplus\mathbb{Z}_2$ & $\mathbb{Z}_2\oplus\mathbb{Z}_2\oplus\mathbb{Z}_4$ & $\mathbb{Z}_2\oplus\mathbb{Z}_2\oplus\mathbb{Z}_2\oplus\mathbb{Z}_2$\\\hline
\end{tabular}
}
\end{center}
depending on whether $m_1,m_2>0$ or $m_1=0$ or $m_2=0$ respectively.

\item If $r\equiv 2,6 \,\,({\rm mod}\,\, 8)$, $N=d_rm$ and
 \[
 \pi_1(\mathcal{N}_{SO(n)}^0(\widehat{Spin(r)}))=
\left\{\begin{array}{ll}
\mathbb{Z}, & \mbox{if $(m,4)=1$,}\\
\mathbb{Z}\times\mathbb{Z}_2, & \mbox{if $(m,4)=2$,}\\
\mathbb{Z}\times\mathbb{Z}_4, & \mbox{if $(m,4)=4$.}
\end{array}
\right.
\]

\item If $r\equiv 3,5 \,\,({\rm mod}\,\, 8)$, $N=d_rm$ and
 \[\pi_1(\mathcal{N}_{SO(N)}^0(\widehat{Spin(r)}))=\mathbb{Z}_2 . \]

 \item If $r\equiv 4 \,\,({\rm mod}\,\, 8)$, $N=d_r(m_1+m_2)$ and
 \begin{eqnarray*}
\pi_1(\mathcal{N}_{SO(N)}^0(\widehat{Spin(r)}))&\cong&\mathbb{Z}_2\oplus \mathbb{Z}_2,\quad \mbox{if $m_1>0$, $m_2>0$,}\\  
\pi_1(\mathcal{N}_{SO(N)}^0(\widehat{Spin(r)^-}))&\cong&
\left\{
\begin{array}{ll}
\mathbb{Z}_2\oplus \mathbb{Z}_2, & \mbox{if $m_1>0$, $m_2=0$, $r>4$,}\\
\mathbb{Z}_2, & \mbox{if $m_1>0$, $m_2=0$, $r=4$,}
\end{array}
\right.\\ 
\pi_1(\mathcal{N}_{SO(N)}^0(\widehat{Spin(r)^+}))&\cong& \left\{
\begin{array}{ll}
\mathbb{Z}_2\oplus \mathbb{Z}_2, & \mbox{if $m_1=0$, $m_2>0$, $r>4$,}\\
\mathbb{Z}_2 & \mbox{if $m_1=0$, $m_2>0$, $r=4$.}
\end{array}
\right.
 \end{eqnarray*}

\end{itemize}
\end{theo}

\qd

\section{Lifting maps to the Spin group}\label{sec: lifts}

In this section, we will check how the generators of the fundamental groups $\pi_1(\mathcal{N}^0_{SO(N)}(S))$ 
map into $\pi_1(SO(N))$.

\begin{theo} \label{theo: lifts}
Let $r\geq 3$. There exist lifts 
\[\begin{array}{ccc}
 &  & Spin(N)\\
 & \nearrow & \downarrow\\
\mathcal{N}^0_{SO(N)}(S) & \longrightarrow & SO(N)
  \end{array}
\]
where $S$ denotes the homomorphic image of $Spin(r)$ in $SO(N)$ 
(either $\widehat{Spin(r)}$ or $\widehat{Spin(r)^\pm}$),
in the following cases:
\begin{itemize} 
 \item \underline{$r\equiv 1,7 \,\,({\rm mod}\,\, 8)$}. 
 \begin{itemize}
  \item For all $m\in \mathbb{N}$.
 \end{itemize}
 \item \underline{$r\equiv 0 \,\,({\rm mod}\,\, 8)$} 
 \begin{itemize}
  \item  For all $m_1,m_2\in\mathbb{N}$ if $r>8$.
    \item For $m_1\equiv m_2\equiv 0$ {\rm (mod 2)} if $r=8$.
 \end{itemize}
 \item \underline{$r\equiv 2,6 \,\,({\rm mod}\,\, 8)$}
    \begin{itemize}
    \item For all $m\in \mathbb{N}$ if $r>6$.
    \item For $m$ even if $r=6$. 
    \end{itemize} 
 \item \underline{$r\equiv 3,5 \,\,({\rm mod}\,\, 8)$} 
    \begin{itemize}
    \item For all $m\in \mathbb{N}$ if $r>3$.
    \item For $m$ even if $r=3$.
    \end{itemize} 
 \item \underline{$r\equiv 4 \,\,({\rm mod}\,\, 8)$}
    \begin{itemize}
    \item For all $m_1,m_2\in \mathbb{N}$ if $r>4$.
    \item For $m_1\equiv m_2\equiv 0$ {\rm (mod 2)} if $r=4$.
    \end{itemize}
\end{itemize}
\end{theo}

The rest of this section is devoted to prove Theorem \ref{theo: lifts} in a case by case analysis.

\subsection{$r\equiv 1,7 \,\,({\rm mod}\,\, 8)$}

Recall 
 \begin{eqnarray*}
\pi_1(\widehat{SO(m)}\widehat{Spin(r)})
 &=&\left\{
 \begin{array}{ll}
\left<(-1,1)\right>\times\left<({\rm vol}_m,-1)\right>\subset Spin(m)\times Spin(r), & \mbox{\rm if $m\geq 4, m\equiv 0$ (mod 4),} \\
\left<({\rm vol}_m,-1)\right>\subset Spin(m)\times Spin(r), & \mbox{\rm if $m\geq 4, m\equiv 2$ (mod 4),}\\
\left<(-1,1)\right>\subset Spin(m)\times Spin(r), & \mbox{\rm if $m\geq 3, m$ is odd,} \\
\{1\}\subset \{1\}\times Spin(r), & \mbox{\rm if $m= 1$,}\\ 
\left<(\pi,-1)\right>\subset \mathbb{R}\times Spin(r), & \mbox{\rm if $m=2$.} 
 \end{array}
\right.\\
 \end{eqnarray*}
Thus, we only need to check the loops in $SO(d_rm)$ which are images of paths joining $(1,1)$ to either 
$(-1,1)$ or $({\rm vol}_m,-1)$ in $Spin(m)\times Spin(r)$ or joining $(0,1)$ to $(\pi,-1)$ in $\mathbb{R}\times Spin(r)$.

\begin{itemize}
 
\item Consider the path
 \begin{eqnarray*}
  \delta_1 : [0,1] &\longrightarrow & Spin(m)\times Spin(r)\\
  t &\mapsto& (\cos(\pi t)+ \sin(\pi t)v_1v_2,1)
 \end{eqnarray*}
joining  $(1,1)$ to $(-1,1)$ in $Spin(m)\times Spin(r)$. It projects to the loop 
\begin{eqnarray*}
  \hat\delta_1 : [0,1] &\longrightarrow & \widehat{SO(m)}\widehat{Spin(r)}\subset SO(N)\\
t&\mapsto& 
\left(
\begin{array}{ccccc}
\cos(2\pi t) & -\sin(2\pi t) &  &  & \\
\sin(2\pi t) & \cos(2\pi t) &  &  & \\
 &  & 1 &  & \\
 &  &  & \ddots & \\
 &  &  &  & 1
\end{array}
\right)_{m\times m}\otimes {\rm Id}_{\Delta_r}, 
\end{eqnarray*}
which contains $2^{[{r\over 2}]}$ blocks
\[\left(
\begin{array}{cc}
\cos(2\pi t) & -\sin(2\pi t)  \\
\sin(2\pi t) & \cos(2\pi t) 
\end{array}
\right).\]
Thus, $\hat\delta_1$ represents $2^{[{r\over 2}]}$ times the generator of $\pi_1(SO(d_rm))$.
Since $r\geq 3$ and $r\equiv 1,7 \,\,({\rm mod}\,\, 8)$, $2^{[{r\over 2}]}$ is divisible by 8 and $\hat\delta_1$  is null homotopic.

\item When $m$ is even and $m\geq 4$, also
consider the path
 \begin{eqnarray*}
  \delta_2 : [0,1] &\longrightarrow & Spin(m)\times Spin(r)\\
  t &\mapsto& (\prod_{j=1}^{m\over 2} \cos(\pi t/2)+ \sin(\pi t/2)v_{2j-1}v_{2j} ,\cos(\pi t)+ \sin(\pi t)e_1e_2)
 \end{eqnarray*}
joining $(1,1)$ to $({\rm vol}_m,-1)$  in $Spin(m)\times Spin(r)$. It projects to the loop 
{\tiny
\begin{eqnarray*}
\hat\delta_2 : [0,1] &\longrightarrow& \widehat{SO(m)}\widehat{Spin(r)}\subset SO(N)\\
t&\mapsto&
\left(
\begin{array}{ccccc}
\cos(\pi t) & -\sin(\pi t) &  &&   \\
\sin(\pi t) & \cos(\pi t) &  & &  \\
 &  & \ddots &  & \\
&&&\cos(\pi t) & -\sin(\pi t) \\
&&&\sin(\pi t) & \cos(\pi t)  
\end{array}
\right)_{m\times m}
\otimes
\left(
\begin{array}{ccccc}
e^{\pi i t} &  &  &&   \\
 & e^{-\pi i t} &  & &  \\
 &  & \ddots &  & \\
&&&e^{\pi i t} &  \\
&&&& e^{-\pi i t}  
\end{array}
\right)_{2^{[{r\over 2}]}\times 2^{[{r\over 2}]}}
\end{eqnarray*}
}
which is similar to
\[
\left(
\begin{array}{ccccc}
e^{\pi i t} &  &  &&   \\
 & e^{-\pi i t} &  & &  \\
 &  & \ddots &  & \\
&&&e^{\pi i t} &  \\
&&&& e^{-\pi i t}  
\end{array}
\right)_{m\times m}
\otimes
\left(
\begin{array}{ccccc}
e^{\pi i t} &  &  &&   \\
 & e^{-\pi i t} &  & &  \\
 &  & \ddots &  & \\
&&&e^{\pi i t} &  \\
&&&& e^{-\pi i t}  
\end{array}
\right)_{2^{[{r\over 2}]}\times 2^{[{r\over 2}]}}.
\]
It contains $2^{[{r\over 2}]-1}$ blocks
\[\left(\begin{array}{ccccc|ccccc}
e^{2\pi i t} &  &  &  &  &  &  &  &  & \\
 & 1 &  &  &  &  &  &  &  & \\
 &  & \ddots &  &  &  &  &  &  & \\
 &  &  & e^{2\pi i t} &  &  &  &  &  & \\
 &  &  &  & 1 &  &  &  &  & \\\hline
 &  &  &  &  & 1 &  &  &  & \\
 &  &  &  &  &  & e^{-2\pi i t} &  &  & \\
 &  &  &  &  &  &  & \ddots &  & \\
 &  &  &  &  &  &  &  & 1 & \\
 &  &  &  &  &  &  &  &  & e^{-2\pi i t}
        \end{array}
\right)_{2m\times 2m}\]
i.e. there are $2^{[{r\over 2}]-1}{m\over 2}=2^{[{r\over 2}]-2}m$ copies of the generator of $\pi_1(SO(d_rm))$.
Since $r\geq 3$ and $r\equiv 1,7 \,\,({\rm mod}\,\, 8)$, $2^{[{r\over 2}]-2}$ is divisible by 2 and $\hat\delta_2$ is null homotopic.

\item For $m=2$, consider the path 
\begin{eqnarray*}
\delta_3:[0,1]&\longrightarrow& \mathbb{R}\times Spin(r) \\
t&\mapsto& (\pi t,\cos(\pi t) + \sin(\pi t) e_1e_2) 
\end{eqnarray*}
joining $(0,1)$ to $(1,-1)$ in $\mathbb{R}\times Spin(r)$, which maps to 
\begin{eqnarray*}
\hat\delta_3:[0,1]&\longrightarrow& \widehat{SO(2)}\widehat{Spin(r)} \subset SO(N)\\
t&\mapsto& \left(\begin{array}{cc}
\cos(\pi t) & -\sin(\pi t)\\
\sin(\pi t) & \cos(\pi t)
                       \end{array}
\right)
\otimes
\left(
\begin{array}{ccccc}
e^{\pi i t} &  &  &&   \\
 & e^{-\pi i t} &  & &  \\
 &  & \ddots &  & \\
&&&e^{\pi i t} &  \\
&&&& e^{-\pi i t}  
\end{array}
\right)_{2^{[{r\over 2}]}\times 2^{[{r\over 2}]}}
\\
&\sim& \left(\begin{array}{cc}
e^{2\pi i t} & \\
 & e^{-2\pi i t}
                       \end{array}
\right)\otimes {\rm Id}_{2^{[{r\over 2}]-1}}\oplus{\rm Id}_{2^{[{r\over 2}]-1}}.
\end{eqnarray*}
This loop represents $2^{[{r\over 2}]-1}$ times the generator of $\pi_1(SO(N))$, which is null homotopic since $2^{[{r\over 2}]-1}$ is divisible by $4$.

\end{itemize}

\subsection{$r\equiv 0 \,\,({\rm mod}\,\, 8)$}

Let $r=8k$, $\{v_1,\ldots,v_{m_1}\}$ and $\{v'_1,\ldots,v'_{m_2}\}$ oriented orthonormal bases of $\mathbb{R}^{m_1}$ and $\mathbb{R}^{m_2}$ respectively. 
Recall the fundamental group generators for $m_1,m_2\geq 3$:
\[
\begin{array}{|ll|c|c|c|c|c|} \hline
&\mbox{Cases} & \pi_1(\widehat{SO(m_1)}\widehat{SO(m_2)}\widehat{Spin(r)}) & (-1,1,1) & (1,-1,1) & ({\rm vol}_{m_1},1,-{\rm vol}_r) & (1,{\rm vol}_{m_2},{\rm vol}_r)\\ \hline
\mbox{(a)}&m_1\equiv 1 (2), m_2\equiv 1 (2) & \mathbb{Z}_2\oplus \mathbb{Z}_2 & \checkmark & \checkmark &  & \\ \hline
\mbox{(b)}&m_1\equiv 0 (4), m_2\equiv 1 (2) & \mathbb{Z}_2\oplus \mathbb{Z}_2\oplus \mathbb{Z}_2 & \checkmark & \checkmark & \checkmark & \\
\mbox{(c)}&m_1\equiv 2 (4), m_2\equiv 1 (2) & \mathbb{Z}_2\oplus \mathbb{Z}_4 & & \checkmark & \checkmark & \\ \hline
\mbox{(d)}&m_1\equiv 1 (2), m_2\equiv 0 (4) & \mathbb{Z}_2\oplus \mathbb{Z}_2\oplus \mathbb{Z}_2 & \checkmark & \checkmark &  &\checkmark \\
\mbox{(e)}&m_1\equiv 1 (2), m_2\equiv 2 (4) & \mathbb{Z}_2\oplus \mathbb{Z}_4 & \checkmark &  &  &\checkmark \\ \hline
\mbox{(f)}&m_1\equiv 0 (4), m_2\equiv 0 (4) & \mathbb{Z}_2\oplus \mathbb{Z}_2\oplus \mathbb{Z}_2\oplus \mathbb{Z}_2 & \checkmark & \checkmark & \checkmark &\checkmark \\
\mbox{(g)}&m_1\equiv 0 (4), m_2\equiv 2 (4) & \mathbb{Z}_2\oplus \mathbb{Z}_2\oplus \mathbb{Z}_4 & \checkmark &  & \checkmark &\checkmark \\
\mbox{(h)}&m_1\equiv 2 (4), m_2\equiv 0 (4) & \mathbb{Z}_2\oplus \mathbb{Z}_4\oplus \mathbb{Z}_2 &  & \checkmark & \checkmark &\checkmark \\
\mbox{(i)}&m_1\equiv 2 (4), m_2\equiv 2 (4) & \mathbb{Z}_4\oplus \mathbb{Z}_4 &  &  & \checkmark & \checkmark \\ \hline
\end{array}
\]

 \begin{itemize}
 
\item For the cases (a), (b), (d), (f) and (g)
consider the path 
\begin{eqnarray*}
  \delta_1 : [0,1] &\longrightarrow& Spin(m_1)\times Spin(m_2)\times Spin(r)\\
  t&\mapsto& (\cos(\pi t) + \sin(\pi t) v_1v_2,1,1) 
\end{eqnarray*}
joining $(1,1,1)$ to $(-1,1,1)$ in $Spin(m_1)\times Spin(m_2)\times Spin(r)$ which projects to the loop 
\begin{eqnarray*}
  \hat\delta_1 : [0,1] &\longrightarrow& (\widehat{SO(m_1)}\times \widehat{SO(m_2)})\widehat{Spin(r)}\subset SO(N)\\
t&\mapsto&
\left(
\begin{array}{ccccc}
\cos(2\pi t) & -\sin(2\pi t) &  &  & \\
\sin(2\pi t) & \cos(2\pi t) &   &  & \\
 &  & 1 &    & \\
 &   &  & \ddots & \\
 &    &  &  & 1
\end{array}
\right)_{m_1\times m_1} \otimes {\rm Id}_{\Delta_r^+} 
\oplus
{\rm Id}_{m_2}\otimes {\rm Id}_{\Delta_r^-}.
\end{eqnarray*}
It contains $2^{{r\over 2}-1}$ copies of the generator of $\pi_1(SO(d_r(m_1+m_2)))$, which is homotopically trivial since
$2^{{r\over 2}-1}$ is divisible by 8.

\item For the cases (a), (b), (c), (d), (f) and (h), 
consider the path 
\begin{eqnarray*}
  \delta_2 : [0,1] &\longrightarrow& Spin(m_1)\times Spin(m_2)\times Spin(r)\\
  t&\mapsto& (1,\cos(\pi t) + \sin(\pi t) v_1'v_2',1) 
\end{eqnarray*}
joining $(1,1,1)$ to $(1,-1,1)$ in $Spin(m_1)\times Spin(m_2)\times Spin(r)$, which projects to the loop 
\begin{eqnarray*}
  \hat\delta_2 : [0,1] &\longrightarrow& (\widehat{SO(m_1)}\times \widehat{SO(m_2)})\widehat{Spin(r)}\subset SO(N)\\
t&\mapsto&{\rm Id}_{m_1}\otimes {\rm Id}_{\Delta_r^+}
\oplus
\left(
\begin{array}{ccccc}
\cos(2\pi t) & -\sin(2\pi t) &  &  & \\
\sin(2\pi t) & \cos(2\pi t) &   &  & \\
 &  & 1 &    & \\
 &   &  & \ddots & \\
 &    &  &  & 1
\end{array}
\right)_{m_2\times m_2} \otimes {\rm Id}_{\Delta_r^-} . 
\end{eqnarray*}
It contains $2^{{r\over 2}-1}$ copies of the generator of $\pi_1(SO(d_r(m_1+m_2)))$, and is homotopically trivial since
$2^{{r\over 2}-1}$ is divisible by 8.

\item For the cases (b), (c), (f), (g), (h) and (i), 
consider the path
\begin{eqnarray*}
  \delta_3 : [0,1] &\longrightarrow& Spin(m_1)\times Spin(m_2)\times Spin(r)\\
  t&\mapsto& (1,\prod_{j=1}^{m_2\over 2}\cos(\pi t/2) + \sin(\pi t/2) v_{2j-1}'v_{2j}',\prod_{l=1}^{r\over 2}\cos(\pi t/2) + \sin(\pi t/2) e_{2l-1}e_{2l}) 
\end{eqnarray*}
joining $(1,1,1)$ to $(1,{\rm vol}_{m_2},{\rm vol}_{r})$ in $Spin(m_1)\times Spin(m_2)\times Spin(r)$. 
It projects to the loop 
\begin{eqnarray*}
  \hat\delta_3 : [0,1] &\longrightarrow& (\widehat{SO(m_1)}\times\widehat{SO(m_2)})\widehat{Spin(r)}\subset SO(N)\\
t&\mapsto&
{\rm Id}_{m_1}\otimes P^+(t)
\oplus
\left(
\begin{array}{ccccc}
\cos(\pi t) & -\sin(\pi t) &  &&   \\
\sin(\pi t) & \cos(\pi t) &  & &  \\
   &  &\ddots &&\\
&&&\cos(\pi t) & -\sin(\pi t) \\
&&&\sin(\pi t) & \cos(\pi t)  
\end{array}
\right)_{m_2\times m_2} \otimes P^-(t) 
\\
&\sim& 
{\rm Id}_{m_1}\otimes P^+(t)
\oplus
\left(
\begin{array}{ccccc}
e^{\pi it} &  &  &&   \\
 & e^{-\pi it} &  & &  \\
   &  &\ddots &&\\
&&&e^{\pi it} &  \\
&&& & e^{-\pi it}  
\end{array}
\right)_{m_2\times m_2} \otimes P^-(t) 
.
\end{eqnarray*}
where 
{\tiny
\begin{eqnarray*}
 P^+(t)
   &=&
   {\rm diag}(
     e^{{2\pi (k)it}},
     \underbrace{e^{{2\pi (k-1)it}},\ldots,e^{{2\pi (k-1)it}}}_{{4k\choose 2}\,\,{\rm times}},
     \underbrace{e^{{2\pi (k-2)it}},\ldots,e^{{2\pi (k-2)it}}}_{{4k\choose 4}\,\,{\rm times}},\ldots, e^{{2\pi (-k)it}}),\\
 P^-(t)
   &=&
   {\rm diag}(
     \underbrace{e^{{(2k-1)}{\pi it}},\ldots,e^{{(2k-1)}{\pi it}}}_{{4k\choose 1}\,\,{\rm times}},
     \underbrace{e^{{(2k-3)}{\pi it}},\ldots,e^{{(2k-3)}{\pi it}}}_{{4k\choose 3}\,\,{\rm times}},\ldots,
     \underbrace{e^{{-(2k-1)}{\pi it}},\ldots,e^{{-(2k-1)}{\pi it}}}_{{4k\choose 4k-1}\,\,{\rm times}}).\\
\end{eqnarray*}
}
Thus, $\hat\delta_3$ containts
{\tiny
\begin{eqnarray*}
&&
m_1\left( k  + {4k\choose 2}(k-1) + {4k\choose 4}(k-2) + \cdots + {4k\choose 2k-2}    \right)  
+ {m_2\over 2}\left( {4k\choose 1}k +  {4k\choose 3}(k-1) + \cdots + {4k\choose 4k-1}(-(k-1))\right)
\\
&=& 
m_1{k(2k-1)\over 8k-2)}{4k\choose 2k} + 2^{4k-3}m_2
\end{eqnarray*}
}
copies of the generator of $\pi_1(SO(d_r(m_1+m_2)))$.
This number of copies is always even except when $k=1$ and $m_1$ is odd.

\item For the cases (d), (e), (f), (g), (h) and (i), 
consider the path 
{\footnotesize
\begin{eqnarray*}
\delta_4: [0,1]&\longrightarrow& Spin(m_1)\times Spin(m_2)\times Spin(r)\\ 
           t &\mapsto& \left(
           \prod_{j=1}^{m_1\over 2}(\cos(\pi t/2) + \sin(\pi t/2)v_{2j-1}v_{2j},1,(\cos(\pi t/2) - \sin(\pi t/2)e_{1}e_{2})\prod_{l=2}^{r\over 2}(\cos(\pi t/2) + \sin(\pi t/2)e_{2l-1}e_{2l} )
           \right)
\end{eqnarray*}
}
joining $(1,1,1)$ to $({\rm vol}_{m_1},1,-{\rm vol}_r)$ in $Spin(m_1)\times Spin(m_2)\times Spin(r)$.
It projects to the loop 
{\footnotesize
\begin{eqnarray*}
\hat\delta_4: [0,1]&\longrightarrow& (\widehat{SO(m_1)}\times\widehat{SO(m_2)})\widehat{Spin(r)}\subset SO(N)\\ 
           t &\mapsto& 
          \left(\begin{array}{ccccc}
\cos(\pi t) & -\sin(\pi t) &  &  & \\
\sin(\pi t) & \cos(\pi t) &  &  & \\
 &  & \ddots &  & \\
 &  &  & \cos(\pi t) & -\sin(\pi t)\\
 &  &  & \sin(\pi t) & \cos(\pi t)
                \end{array}
\right)_{m_1\times m_1}\otimes Q^+(t) 
\oplus
{\rm Id}_{m_2}\otimes Q^-(t)\\
&\sim& 
          \left(\begin{array}{ccccc}
e^{\pi i t} &  &  &  & \\
 & e^{-\pi i t} &  &  & \\
 &  & \ddots &  & \\
 &  &  & e^{\pi i t} & \\
 &  &  &  & e^{-\pi i t}
                \end{array}
\right)_{m_1\times m_1}\otimes Q^+(t) 
\oplus
{\rm Id}_{m_2}\otimes Q^-(t)\\
\end{eqnarray*}
}
where 
{\tiny
\begin{eqnarray*}
 Q^+(t)
   &=&
   {\rm diag}(
     \underbrace{e^{{(2k-1)\pi it}},\ldots,e^{{(2k-1)\pi it}}}_{{4k\choose 1}\,\,{\rm times}},
     \underbrace{e^{{(2k-3)\pi it}},\ldots,e^{{(2k-3)\pi it}}}_{{4k\choose 3}\,\,{\rm times}},
     \underbrace{e^{{(2k-5)\pi it}},\ldots,e^{{(2k-5)\pi it}}}_{{4k\choose 5}\,\,{\rm times}},\ldots,
     \underbrace{e^{{-(2k-1)\pi it}},\ldots,e^{{-(2k-1)\pi it}}}_{{4k\choose 4k-1}\,\,{\rm times}}),\\
 Q^-(t)
   &=&
   {\rm diag}(
     e^{{2\pi kit}},
     \underbrace{e^{{2\pi (k-1)it}},\ldots,e^{{2\pi (k-1)it}}}_{{4k\choose 2}\,\,{\rm times}},
     \underbrace{e^{{2\pi (k-2)it}},\ldots,e^{{2\pi (k-2)it}}}_{{4k\choose 4}\,\,{\rm times}},
     \ldots,
     e^{{2\pi (-k)it}}).
\end{eqnarray*}
}
Thus, we have
{\footnotesize
\begin{eqnarray*}
&&{m_1\over 2} \left[ k {4k\choose 1} + (k-1) {4k\choose 3}+\cdots+(-(k-1)){4k\choose 4k-1}\right]
+ m_2 \left[ k {4k\choose 0} + (k-1){4k\choose 2}+\cdots + 1{4k\choose 2k-2}\right]\\
&=&
2^{4k-3}m_1+m_2{k(2k-1)\over 8k-2}{4k\choose 2k},
\end{eqnarray*}
}
which is odd only if $k=1$ and $m_2\equiv 1$ (mod 2).

\end{itemize}

The cases in which either $m_1\leq 2$ or $m_2\leq 2$ are treated similarly.

\subsection{$r\equiv 2,6 \,\,({\rm mod}\,\, 8)$}

Recall
 \[\pi_1(\widehat{U(m)}\widehat{Spin(r)})=
\left\{\begin{array}{ll}
\mathbb{Z}, & \mbox{\rm if $(m,4)=1$,}\\
\mathbb{Z}\oplus\mathbb{Z}_2, & \mbox{\rm if $(m,4)=2$,}\\
\mathbb{Z}\oplus\mathbb{Z}_4, & \mbox{\rm if $(m,4)=4$.}
\end{array}
\right.
\]
In every case, the fundamental group has generators
\[\left( {2\pi\over m}, e^{-{2\pi i \over m}}{\rm Id}_{m}, 1\right),\quad \left({\pi\over 2}, {\rm Id}_{m}, -{\rm vol}\right).\]

\begin{itemize}
\item 
Consider the path 
\begin{eqnarray*}
 \delta_1: [0,1] &\longrightarrow& \mathbb{R}\times SU(m)\times Spin(r) \\
  t &\mapsto& \left({2\pi t\over m},
  \left(\begin{array}{ccccc}
e^{-{2\pi i t\over m}} &  &  &  & \\
 & e^{-{2\pi i t\over m}} &  &  & \\
 &  & \ddots &  & \\
 &  &  & e^{-{2\pi i t\over m}} & \\
 &  &  &  & e^{{2\pi i (m-1)t\over m}}
        \end{array}
  \right)_{m\times m}
  , 1\right)
\end{eqnarray*}
joining $(0, {\rm Id}_{m \times m},1)$ to
$\left( {2\pi\over m}, e^{-{2\pi i \over m}}{\rm Id}_{m}, 1\right)$ in $\mathbb{R}\times SU(m)\times Spin(r)$.
This path gets mapped to the loop 
\begin{eqnarray*}
 \hat\delta_1: [0,1] &\longrightarrow& \widehat{U(m)}\widehat{Spin(r)} \subset SO(N)\\
  t &\mapsto& 
  \left(\begin{array}{c|c}
  \left(\begin{array}{ccccc}
1 &  &  &  & \\
 & 1 &  &  & \\
 &  & \ddots &  & \\
 &  &  & 1 & \\
 &  &  &  & e^{2\pi i t}
        \end{array}
  \right)_{m\times m}
  \otimes  {\rm Id}_{\Delta_r^+} & 0\\\hline
0&    \left(\begin{array}{ccccc}
1 &  &  &  & \\
 & 1 &  &  & \\
 &  & \ddots &  & \\
 &  &  & 1 & \\
 &  &  &  & e^{-2\pi i t}
        \end{array}
  \right)_{m\times m}
  \otimes  {\rm Id}_{\Delta_r^-} 
\end{array}
\right)  
\end{eqnarray*}
if $r\equiv 2 \,\,({\rm mod}\,\, 8)$, and to
\begin{eqnarray*}
 \hat\delta_1: [0,1] &\longrightarrow& \widehat{U(m)}\widehat{Spin(r)} \subset SO(N) \\
  t &\mapsto& 
  \left(\begin{array}{c|c}
  \left(\begin{array}{ccccc}
1 &  &  &  & \\
 & 1 &  &  & \\
 &  & \ddots &  & \\
 &  &  & 1 & \\
 &  &  &  & e^{-2\pi i t}
        \end{array}
  \right)_{m\times m}
  \otimes  {\rm Id}_{\Delta_r^+} & 0\\\hline
0&    \left(\begin{array}{ccccc}
1 &  &  &  & \\
 & 1 &  &  & \\
 &  & \ddots &  & \\
 &  &  & 1 & \\
 &  &  &  & e^{2\pi i t}
        \end{array}
  \right)_{m\times m}
  \otimes  {\rm Id}_{\Delta_r^-} 
\end{array}
\right)  
\end{eqnarray*}
if $r\equiv 6 \,\,({\rm mod}\,\, 8)$.
Since $\dim(\Delta_r^\pm)=2^{{r\over 2}-1}$, we have $2^{{r\over 2}-1}$  blocks of the form
\[\left(\begin{array}{cc}
e^{2\pi i t} & \\
 & e^{-2\pi i t}
        \end{array}
\right) \sim
\left(\begin{array}{cc}
\cos(2\pi  t) & -\sin(2\pi  t)\\
\sin(2\pi  t) & \cos(2\pi  t)
        \end{array}
\right) 
\]
i.e. $2^{{r\over 2}-1}$ times the generator of $\pi_1(SO(d_rm))$. Since 
$r\geq 6$,  $2^{{r\over 2}-1}$ is divisible by 4. Hence, $\hat\delta_1$ is null homotopic.

\item 
Consider the path 
\begin{eqnarray*}
 \delta_2: [0,1] &\longrightarrow& \mathbb{R}\times SU(m)\times Spin(r) \\
  t &\mapsto& \left({\pi t\over 2},
  {\rm Id}_{m}
  , \prod_{j=1}^{{r\over 2}}(\cos(\pi t/2)-\sin(\pi t/2)e_{2j-1}e_{2j})\right)
\end{eqnarray*}
joining $(0, {\rm Id}_{m \times m},1)$ to
$\left({\pi\over 2}, {\rm Id}_{m}, -{\rm vol}\right)$ in $\mathbb{R}\times SU(m)\times Spin(r)$.
This path gets mapped to the loop 
\begin{eqnarray*}
 \hat\delta_2: [0,1] &\longrightarrow& \widehat{U(m)}\widehat{Spin(r)} \subset SO(N)\\
  t &\mapsto& 
  \left(\begin{array}{c|c}
         e^{\pi i t\over 2}{\rm Id}_{m}
  \otimes  P^+(t) & 0\\\hline
0&   e^{-{\pi i t\over 2}}{\rm Id}_{m} \otimes  P^-(t)
\end{array}
\right)  
\end{eqnarray*}
if $r\equiv 2 \,\,({\rm mod}\,\, 8)$, and
\begin{eqnarray*}
 \hat\delta_2: [0,1] &\longrightarrow& \widehat{U(m)}\widehat{Spin(r)}\subset SO(N) \\
  t &\mapsto& 
  \left(\begin{array}{c|c}
         e^{-{\pi i t\over 2}}{\rm Id}_{m}
  \otimes  P^+(t) & 0\\\hline
0&   e^{{\pi i t\over 2}}{\rm Id}_{m} \otimes  P^-(t)
\end{array}
\right)  
\end{eqnarray*}
if $r\equiv 6 \,\,({\rm mod}\,\, 8)$, 
where 
\begin{eqnarray*}
 P^+(t)
   &=&
   {\rm diag}(
     e^{-{r\over 2}{\pi it\over 2}},
     \underbrace{e^{-{({r\over 2}-4)}{\pi it\over 2}},\ldots,e^{-{({r\over 2}-4)}{\pi it\over 2}}}_{{r/2\choose 2}\,\,{\rm times}},
     \underbrace{e^{-{({r\over 2}-8)}{\pi it\over 2}},\ldots,e^{-{({r\over 2}-8)}{\pi it\over 2}}}_{{r/2\choose 4}\,\,{\rm times}},\ldots),\\
 P^-(t)
   &=&
   {\rm diag}(
     \underbrace{e^{-{({r\over 2}-2)}{\pi it\over 2}},\ldots,e^{-{({r\over 2}-2)}{\pi it\over 2}}}_{{r/2\choose 1}\,\,{\rm times}},
     \underbrace{e^{-{({r\over 2}-6)}{\pi it\over 2}},\ldots,e^{-{({r\over 2}-6)}{\pi it\over 2}}}_{{r/2\choose 3}\,\,{\rm times}},
     \underbrace{e^{-{({r\over 2}-10)}{\pi it\over 2}},\ldots,e^{-{({r\over 2}-10)}{\pi it\over 2}}}_{{r/2\choose 5}\,\,{\rm times}},\ldots).
\end{eqnarray*}

\begin{itemize}
 \item If $r\equiv 2 \,\,({\rm mod}\,\, 8)$, $r=8k+2$ with $k\geq 1$. Then ${r\over 2}= 4k+1$.
We have 
{\footnotesize
\[m\left[k{4k+1\choose0} + (k-1){4k+1\choose2} + (k-2){4k+1\choose4} + \cdots + (-k+1){4k+1\choose4k-2}  + (-k){4k+1\choose4k}\right]  = -m 2^{4k-2}.\]
}copies of the generator of $\pi_1(SO(d_rm))$, which is even and $\hat\delta_1$ is null homotopic.

 \item If $r\equiv 6 \,\,({\rm mod}\,\, 8)$, $r=8k+6$ with $k\geq 1$. Then ${r\over 2}= 4k+3$.
Thus, we have 
{\footnotesize
\[m\left[(k+1){4k+3\choose0} + (k){4k+3\choose2} + (k-1){4k+3\choose4} + \cdots + (-k+1){4k+3\choose4k}  + (-k){4k+3\choose4k+2} \right] = m 2^{4k}\]
}copies of the generator of $\pi_1(SO(d_rm))$.
If $k\geq 1$, this number is always
even and $\hat\delta_2$ is null homotopic.
On the other hand, if $r=6$ ($k=0$), then the parity of the number depends on $m$.

\end{itemize}

\end{itemize}

\subsection{$r\equiv 3,5 \,\,({\rm mod}\,\, 8)$}

Recall
\[\pi_1(\widehat{Sp(m)}\widehat{Spin(r)})=\mathbb{Z}_2 =\left<(-{\rm Id}_{2m},-1)\right> . \]
Thus, consider the path
\begin{eqnarray*}
\delta: [0,1]&\longrightarrow& Sp(m)\times Spin(r)\\ 
           t &\mapsto& \left(
           \left(
\begin{array}{lllll}
e^{\pi it} &  &  &  & \\
 & e^{-\pi it} &  &  & \\
 &  & \ddots &  & \\
 &  &  & e^{\pi it} & \\
 &  &  &  & e^{-\pi it}
\end{array}
\right)_{2m\times 2m}, \cos(\pi t) + \sin(\pi t)e_1e_2           
           \right)
\end{eqnarray*}
joining $({\rm Id}_{2m},1)$ to $(-{\rm Id}_{2m},-1)$ in $Sp(m)\times Spin(r)$.
It projects to the loop in $\widehat{Sp(m)}\widehat{Spin(r)}\subset SO(d_rm)$
\[           \left(
\begin{array}{ccccc}
e^{\pi it} &  &  &  & \\
 & e^{-\pi it} &  &  & \\
 &  & \ddots &  & \\
 &  &  & e^{\pi it} & \\
 &  &  &  & e^{-\pi it}
\end{array}
\right)_{2m\times 2m}\otimes \left(
\begin{array}{lllll}
e^{\pi it} &  &  &  & \\
 & e^{-\pi it} &  &  & \\
 &  & \ddots &  & \\
 &  &  & e^{\pi it} & \\
 &  &  &  & e^{-\pi it}
\end{array}
\right)_{2^{[{r\over 2}]}\times 2^{[{r\over 2}]}}          . 
\]
It has
$2^{[{r\over 2}]-1}m$ blocks of the form
\[\left(
\begin{array}{cc}
e^{2\pi it} & \\
 & e^{-2\pi it}
\end{array}
\right)\sim 
\left(\begin{array}{ll}
\cos(2\pi t) & -\sin(2\pi t)\\
\sin(2\pi t) & \cos(2\pi t)
      \end{array}
\right),\]
i.e.
$2^{[{r\over 2}]-1}m$ times the generator of $\pi_1(SO(d_rm))=\mathbb{Z}_2$.
Note that $2^{[{r\over 2}]-1}m$ is divisible by 2 if $r>3$.

\subsection{$r\equiv 4 \,\,({\rm mod}\,\, 8)$}

Let $r=8k+4$.
Recall
{\footnotesize
 \begin{eqnarray*}
\pi_1((\widehat{Sp(m_1)}\times\widehat{Sp(m_2)})\widehat{Spin(r)})
   &=&
  \left<(-{\rm Id}_{2m_1}, -{\rm Id}_{2m_2},-1)\right> \times
     \left<({\rm Id}_{2m_1},-{\rm Id}_{2m_2},{\rm vol}_r)\right>\subset Sp(m_1)\times Sp(m_2)\times Spin(r). \\
\pi_1(\widehat{Sp(m_1)}\widehat{Spin(r)^-})
&\cong&
\left\{
\begin{array}{ll}
\left<(-{\rm Id}_{2m_1},-{\rm vol}_r),(-{\rm Id}_{2m_1},-1)\right>\subset Sp(m_1)\times Spin(r), & \mbox{if $m_1>0$, $m_2=0$, $r>4$,}\\
\left<(-{\rm Id}_{2m_1},-{\rm vol}_4)\right>\subset Sp(m_1)\times Spin(4), & \mbox{if $m_1>0$, $m_2=0$, $r=4$,}
\end{array}
\right.\\ 
\pi_1(\widehat{Sp(m_2)}\widehat{Spin(r)^+})
&\cong&
\left\{
\begin{array}{ll}
\left<(-{\rm Id}_{2m_2},{\rm vol}_r),(-{\rm Id}_{2m_2},-1)\right>\subset Sp(m_2)\times Spin(r), & \mbox{if $m_1=0$, $m_2>0$, $r>4$,}\\
\left<(-{\rm Id}_{2m_2},{\rm vol}_4)\right>\subset Sp(m_2)\times Spin(4), & \mbox{if $m_1=0$, $m_2>0$, $r=4$,}
\end{array}
\right.\\ 
 \end{eqnarray*}
 }

\begin{itemize}
 \item 
Consider the  path 
{\tiny
\begin{eqnarray*}
\delta_1: [0,1]&\longrightarrow& Sp(m_1)\times Sp(m_2)\times Spin(r)\\ 
           t &\mapsto& \left(
           \left(
\begin{array}{ccccc}
e^{\pi it} &  &  &  & \\
 & e^{-\pi it} &  &  & \\
 &  & \ddots &  & \\
 &  &  & e^{\pi it} & \\
 &  &  &  & e^{-\pi it}
\end{array}
\right)_{2m_1\times 2m_1}, 
           \left(
\begin{array}{ccccc}
e^{\pi it} &  &  &  & \\
 & e^{-\pi it} &  &  & \\
 &  & \ddots &  & \\
 &  &  & e^{\pi it} & \\
 &  &  &  & e^{-\pi it}
\end{array}
\right)_{2m_2\times 2m_2}, 
\cos(\pi t) + \sin(\pi t)e_1e_2           
           \right)
\end{eqnarray*}
}
joining $({\rm Id}_{2m_1},{\rm Id}_{2m_2},1)$ and $(-{\rm Id}_{2m_1}, -{\rm Id}_{2m_2},-1)$ in $Sp(m_1)\times Sp(m_2)\times Spin(r)$.
It maps to the loop 
\begin{eqnarray*} 
\hat\delta_1 : [0,1] &\longrightarrow& (\widehat{Sp(m_1)}\times\widehat{Sp(m_2)})\widehat{Spin(r)}\subset SO(N)\\ 
t&\mapsto&{\rm diag}(e^{\pi i t},e^{-\pi i t},\ldots,e^{\pi i t},e^{-\pi i t})_{2m_1\times 2m_1}\otimes 
{\rm diag}(e^{\pi i t},e^{-\pi i t},\ldots,e^{\pi i t},e^{-\pi i t})_{2^{{r\over 2}-1}\times 2^{{r\over 2}-1}} \\
&&\oplus {\rm diag}(e^{\pi i t},e^{-\pi i t},\ldots,e^{\pi i t},e^{-\pi i t})_{2m_2\times 2m_2}\otimes 
{\rm diag}(e^{\pi i t},e^{-\pi i t},\ldots,e^{\pi i t},e^{-\pi i t})_{2^{{r\over 2}-1}\times 2^{{r\over 2}-1}} \\
&\sim&2^{{r\over 2}-2}{\rm diag}(e^{2\pi i t},e^{-2\pi i t},\ldots,e^{2\pi i t},e^{-2\pi i t})_{2m_1\times 2m_1}
\oplus 2^{{r\over 2}-2} {\rm Id}_{2m_1} \\
&&\oplus 2^{{r\over 2}-2}{\rm diag}(e^{2\pi i t},e^{-2\pi i t},\ldots,e^{2\pi i t},e^{-2\pi i t})_{2m_2\times 2m_2}
\oplus 2^{{r\over 2}-2} {\rm Id}_{2m_2} 
\end{eqnarray*}
which is homotopically equivalent to 
$(m_1+m_2)2^{{r\over 2}-2}$
times the generator of $\pi_1(SO(d_r(m_1+m_2)))$. 
Hence, $\hat\delta_1$ is null homotopic if either $r\not =4$ or $r=4$ and $m_1+ m_2\equiv 0$ (mod 2).

\item 
Consider the path 
{\tiny
\begin{eqnarray*}
\delta_2: [0,1]&\longrightarrow& Sp(m_1)\times Sp(m_2)\times Spin(r)\\ 
           t &\mapsto& \left(
           {\rm Id}_{2m_1},
           \left(
\begin{array}{ccccc}
e^{\pi it} &  &  &  & \\
 & e^{-\pi it} &  &  & \\
 &  & \ddots &  & \\
 &  &  & e^{\pi it} & \\
 &  &  &  & e^{-\pi it}
\end{array}
\right)_{2m_2\times 2m_2}, 
\prod_{j=1}^{r\over 2}\cos(\pi t/2) + \sin(\pi t/2)e_{2j-1}e_{2j}           
           \right).
\end{eqnarray*}
}
joining $({\rm Id}_{2m_1},{\rm Id}_{2m_2},1)$ to
$({\rm Id}_{2m_1},-{\rm Id}_{2m_2},{\rm vol}_r)$.
It maps to the loop 
\begin{eqnarray*}
\hat\delta_2 &\longrightarrow& (\widehat{Sp(m_1)}\times \widehat{Sp(m_2)})\widehat{Spin(r)}\subset SO(N)\\
t&\mapsto&{\rm Id}_{2m_1}\otimes P^-(t)\oplus {\rm diag}(e^{\pi i t},e^{-\pi i t},\ldots,e^{\pi i t},e^{-\pi i t})_{2m_2\times 2m_2}\otimes 
P^+(t)  ,
\end{eqnarray*}
where 
\begin{eqnarray*}
 P^+(t)
   &=&
   {\rm diag}(
     e^{{(2k+1)\pi it}},
     \underbrace{e^{{(2k-1)\pi it}},\ldots,e^{{(2k-1)\pi it}}}_{{4k+2\choose 2}\,\,{\rm times}},
     \underbrace{e^{{(2k-3)\pi it}},\ldots,e^{{(2k-3)\pi it}}}_{{4k+2\choose 4}\,\,{\rm times}},\ldots,e^{{-(2k+1)\pi it}}),\\
 P^-(t)
   &=&
   {\rm diag}(
     \underbrace{e^{{2\pi kit}},\ldots,e^{{2\pi kit}}}_{{4k+2\choose 1}\,\,{\rm times}},
     \underbrace{e^{{2\pi (k-1)it}},\ldots,e^{{2\pi (k-1)it}}}_{{4k+2\choose 3}\,\,{\rm times}},\ldots,
     \underbrace{e^{{2\pi (-k)it}},\ldots,e^{{2\pi (-k)it}}}_{{4k+2\choose 4k+1}\,\,{\rm times}}).
\end{eqnarray*}
This loop is homotopically equivalent (mod 2) to 
{\tiny
\begin{eqnarray*}
&&m_2\left({4k+2\choose 0}(k+1) + {4k+2\choose 2} k + {4k+2\choose 4} (k-1) + \cdots + {4k+2\choose 4k}(-(k-1)) + {4k+2\choose 4k+2}(-k)  \right)
+2m_1\left( {4k+2\choose 1}k + {4k+2\choose 3}(k-1) + \cdots + {4k+2\choose 2k-1}    \right)\\
&=&m_216^k
+2m_1{k(2k+1)\over 8k+2}{4k+2\choose 2k+1} 
\end{eqnarray*}
}
times the generator of $\pi_1(SO(N))$, which is trivial if $k\geq 1$.

 \item 
Consider the  path 
{\tiny
\begin{eqnarray*}
\delta_3: [0,1]&\longrightarrow& Sp(m_1)\times Spin(r)\\ 
           t &\mapsto& \left(
           \left(
\begin{array}{ccccc}
e^{\pi it} &  &  &  & \\
 & e^{-\pi it} &  &  & \\
 &  & \ddots &  & \\
 &  &  & e^{\pi it} & \\
 &  &  &  & e^{-\pi it}
\end{array}
\right)_{2m_1\times 2m_1}, 
\cos(\pi t) + \sin(\pi t)e_1e_2           
           \right)
\end{eqnarray*}
}
joining $({\rm Id}_{2m_1},1)$ and $(-{\rm Id}_{2m_1}, -1)$ in $Sp(m_1)\times Spin(r)$.
It maps to the loop 
\begin{eqnarray*} 
\hat\delta_3 : [0,1] &\longrightarrow& \widehat{Sp(m_1)}\widehat{Spin(r)^-}\subset SO(N)\\ 
t&\mapsto&{\rm diag}(e^{\pi i t},e^{-\pi i t},\ldots,e^{\pi i t},e^{-\pi i t})_{2m_1\times 2m_1}\otimes 
{\rm diag}(e^{\pi i t},e^{-\pi i t},\ldots,e^{\pi i t},e^{-\pi i t})_{2^{{r\over 2}-1}\times 2^{{r\over 2}-1}} \\
&\sim&2^{{r\over 2}-2}{\rm diag}(e^{2\pi i t},e^{-2\pi i t},\ldots,e^{2\pi i t},e^{-2\pi i t})_{2m_1\times 2m_1}
\oplus 2^{{r\over 2}-2} {\rm Id}_{2m_1} \\
\end{eqnarray*}
which is homotopically equivalent to 
$m_12^{{r\over 2}-2}$
times the generator of $\pi_1(SO(d_r(m_1)))$. 
Hence, $\hat\delta_3$ is null homotopic since $r\not = 4$.

 \item 
Consider the  path 
{\tiny
\begin{eqnarray*}
\delta_4: [0,1]&\longrightarrow&  Sp(m_2)\times Spin(r)\\ 
           t &\mapsto& \left(
           \left(
\begin{array}{ccccc}
e^{\pi it} &  &  &  & \\
 & e^{-\pi it} &  &  & \\
 &  & \ddots &  & \\
 &  &  & e^{\pi it} & \\
 &  &  &  & e^{-\pi it}
\end{array}
\right)_{2m_2\times 2m_2}, 
\cos(\pi t) + \sin(\pi t)e_1e_2           
           \right)
\end{eqnarray*}
}
joining $({\rm Id}_{2m_2},1)$ and $( -{\rm Id}_{2m_2},-1)$ in $ Sp(m_2)\times Spin(r)$.
It maps to the loop 
\begin{eqnarray*} 
\hat\delta_4 : [0,1] &\longrightarrow& \widehat{Sp(m_2)}\widehat{Spin(r)^+}\subset SO(N)\\ 
t&\mapsto&
{\rm diag}(e^{\pi i t},e^{-\pi i t},\ldots,e^{\pi i t},e^{-\pi i t})_{2m_2\times 2m_2}\otimes 
{\rm diag}(e^{\pi i t},e^{-\pi i t},\ldots,e^{\pi i t},e^{-\pi i t})_{2^{{r\over 2}-1}\times 2^{{r\over 2}-1}} \\
&\sim&
2^{{r\over 2}-2}{\rm diag}(e^{2\pi i t},e^{-2\pi i t},\ldots,e^{2\pi i t},e^{-2\pi i t})_{2m_2\times 2m_2}
\oplus 2^{{r\over 2}-2} {\rm Id}_{2m_2} 
\end{eqnarray*}
which is homotopically equivalent to 
$m_22^{{r\over 2}-2}$
times the generator of $\pi_1(SO(N))$. 
Hence, $\hat\delta_4$ is null homotopic since $r\not =4$.

\item 
Consider the path 
{\tiny
\begin{eqnarray*}
\delta_5: [0,1]&\longrightarrow& Sp(m_2)\times Spin(r)\\ 
           t &\mapsto& \left(
           \left(
\begin{array}{ccccc}
e^{\pi it} &  &  &  & \\
 & e^{-\pi it} &  &  & \\
 &  & \ddots &  & \\
 &  &  & e^{\pi it} & \\
 &  &  &  & e^{-\pi it}
\end{array}
\right)_{2m_2\times 2m_2}, 
\prod_{j=1}^{r\over 2}\cos(\pi t/2) + \sin(\pi t/2)e_{2j-1}e_{2j}           
           \right).
\end{eqnarray*}
}
joining $({\rm Id}_{2m_2},1)$ to
$(-{\rm Id}_{2m_2},{\rm vol}_r)$.
It maps to the loop 
\begin{eqnarray*}
\hat\delta_5 &\longrightarrow& \widehat{Sp(m_2)}\widehat{Spin(r)^+}\subset SO(N)\\
t&\mapsto& {\rm diag}(e^{\pi i t},e^{-\pi i t},\ldots,e^{\pi i t},e^{-\pi i t})_{2m_2\times 2m_2}\otimes 
P^+(t)  ,
\end{eqnarray*}
where 
\begin{eqnarray*}
 P^+(t)
   &=&
   {\rm diag}(
     e^{{(2k+1)\pi it}},
     \underbrace{e^{{(2k-1)\pi it}},\ldots,e^{{(2k-1)\pi it}}}_{{4k+2\choose 2}\,\,{\rm times}},
     \underbrace{e^{{(2k-3)\pi it}},\ldots,e^{{(2k-3)\pi it}}}_{{4k+2\choose 4}\,\,{\rm times}},\ldots,e^{{-(2k+1)\pi it}}),
\end{eqnarray*}
This loop is homotopically equivalent (mod 2) to 
{\tiny
\begin{eqnarray*}
m_2\left({4k+2\choose 0}(k+1) + {4k+2\choose 2} k + {4k+2\choose 4} (k-1) + \cdots + {4k+2\choose 4k}(-(k-1)) + {4k+2\choose 4k+2}(-k)  \right)
=m_216^k
\end{eqnarray*}
}
times the generator of $\pi_1(SO(N))$, which is trivial if $k\geq 1$.

\item 
Consider the path 
{\tiny
\begin{eqnarray*}
\delta_6: [0,1]&\longrightarrow& Sp(m_1)\times Spin(r)\\ 
           t &\mapsto& \left(
           \left(
\begin{array}{ccccc}
e^{\pi it} &  &  &  & \\
 & e^{-\pi it} &  &  & \\
 &  & \ddots &  & \\
 &  &  & e^{\pi it} & \\
 &  &  &  & e^{-\pi it}
\end{array}
\right)_{2m_1\times 2m_1}, 
(\cos(\pi t/2) - \sin(\pi t/2)e_{1}e_{2})\prod_{j=2}^{r\over 2}\cos(\pi t/2) + \sin(\pi t/2)e_{2j-1}e_{2j}           
           \right).
\end{eqnarray*}
}
joining $({\rm Id}_{2m_1},1)$ to
$(-{\rm Id}_{2m_1},-{\rm vol}_r)$.
It maps to the loop 
\begin{eqnarray*}
\hat\delta_6 &\longrightarrow& \widehat{Sp(m_2)}\widehat{Spin(r)^-}\subset SO(N)\\
t&\mapsto& {\rm diag}(e^{\pi i t},e^{-\pi i t},\ldots,e^{\pi i t},e^{-\pi i t})_{2m_2\times 2m_2}\otimes 
Q^-(t)  ,
\end{eqnarray*}
where  
\begin{eqnarray*}
 Q^-(t)
   &=&
   {\rm diag}(
     e^{{(2k+1)\pi it}},
     \underbrace{e^{{(2k-1)\pi it}},\ldots,e^{{(2k-1)\pi it}}}_{{4k+2\choose 2}\,\,{\rm times}},
     \underbrace{e^{{(2k-3)\pi it}},\ldots,e^{{(2k-3)\pi it}}}_{{4k+2\choose 4}\,\,{\rm times}},\ldots,e^{{-(2k+1)\pi it}}),
\end{eqnarray*}
This loop is homotopically equivalent (mod 2) to 
{\tiny
\begin{eqnarray*}
m_1\left({4k+2\choose 0}(k+1) + {4k+2\choose 2} k + {4k+2\choose 4} (k-1) + \cdots + {4k+2\choose 4k}(-(k-1)) + {4k+2\choose 4k+2}(-k)  \right)
=m_116^k
\end{eqnarray*}
}
times the generator of $\pi_1(SO(N))$, which is trivial if $k\geq 1$.

\end{itemize}

\qd

{\small
\renewcommand{\baselinestretch}{0.5}
\newcommand{\bi}{\vspace{-.05in}\bibitem} }

\end{document}